\documentclass[reqno]{amsart}
\usepackage{amsfonts,amssymb,amsthm,amsmath,mathrsfs}

\newtheorem{theorem}{Theorem}[section]
\newtheorem{lemma}[theorem]{Lemma}
\newtheorem{proposition}[theorem]{Proposition}
\newtheorem{corollary}[theorem]{Corollary}

\theoremstyle{definition}
\newtheorem{defn}[theorem]{Definition}
\newtheorem{example}[theorem]{Example}
\newtheorem{remark}[theorem]{Remark}

\newcommand{\Case}[2]    {{\em Case~#1:~#2\/}}
\newcommand{\prfpart}[1] {{\rm (#1)} }

\newcommand{\Aa}{{\mathbb A}}
\newcommand{\Nn}{{\mathbb N}}
\newcommand{\Pp}{{\mathbb P}}

\newcommand{\PV}{{\mathcal V}}               
\newcommand{\PS}{{\mathcal X}}               

\newcommand{\SVA}{\tilde{\mathcal S}}        
\newcommand{\PVA}{\tilde{\mathcal V}}

\newcommand{\SVAP}{\SVA'}                       
\newcommand{\STU}{S_{T,U}}                      

\newcommand{\A}         {{\mathcal A}}          
\newcommand{\bs}[1]     {\big\{#1\big\}}        
\newcommand{\cmp}[1]    {\overline{#1}}         
\newcommand{\defterm}[1]{{\it #1\/}}
\newcommand{\EE}        {{\bf E}}
\newcommand{\edg}       {\EE}
\newcommand{\fld}       {{\bf k}}
\newcommand{\graft}     {\:\nearrow\:}
\newcommand{\graftat}[1]{\,{\underset{#1}{\nearrow}}\;}

\newcommand{\dju}       {\sqcup}

\newcommand{\isom}      {\cong}
\newcommand{\LL}        {{\mathcal L}}          
\newcommand{\Pic}       {{\bf P}}
\newcommand{\mm}        {{\bf m}}
\newcommand{\sd}        {\#}                    

\newcommand{\sm}        {\setminus}
\newcommand{\SR}[1]     {\Delta({#1})}          
\newcommand{\st}        {~\mid~}                
\newcommand{\surj}      {\twoheadrightarrow}
\newcommand{\TP}        {{\tau}}                
\newcommand{\TT}        {{\bf T}}
\newcommand{\x}         {\times}

\newcommand{\Adm}       {\mathop{{\rm Adm}}\nolimits}
\newcommand{\Aug}       {\mathop{{\rm Aug}}\nolimits}
\newcommand{\BT}        {\mathop{{\rm BPP}}\nolimits}           
\newcommand{\Ch}        {\mathop{{\rm Ch}}\nolimits}            
\newcommand{\codim}     {\mathop{{\rm codim}}\nolimits}
\newcommand{\DPT}       {\mathop{{\rm DPT}}\nolimits}           
\newcommand{\dpt}       {d}             

\newcommand{\ini}       {\mathop{{\rm in}}\nolimits}
\newcommand{\Match}     {M}
\newcommand{\Rd}        {\mathop{{\rm Rd}}\nolimits}
\newcommand{\rr}        {\mathop{{\rm rr}}\nolimits}            
\newcommand{\rt}[1]     {{\mathop{\rm rt}\nolimits}(#1)}
\newcommand{\shs}       {\mathop{\rm sh}\nolimits}
\newcommand{\SST}       {\mathop{\rm Cpl}\nolimits}   
\newcommand{\Spec}      {\mathop{\rm Spec}\nolimits}
\newcommand{\tmin}      {\mathop{{\rm min}_2}}

\newcommand{\trav}      {\mathop{\rm trav}\nolimits}
\newcommand{\val}       {\mathop{\rm val}\nolimits}          

\newcommand{\maxv}{\max}
\newcommand{\minv}{\min}
\newcommand{\gtv}{>}
\newcommand{\ltv}{<}
\newcommand{\llex}      {<}             
\newcommand{\glex}      {>}             

\newcommand{\grlex}     {\mathop{>}\limits_{\rm rlex}}
\newcommand{\gint}      {>}             
\newcommand{\lfac}      {\prec}         
\newcommand{\gfac}      {\succ}         

\newcommand{\putvertex}[1]{\put(#1){\circle*{5}}}
\newcommand{\puttext}[2]{\put(#1){\makebox(0,0){#2}}}

\begin{document}

\title{The Slopes Determined by $n$ Points in the Plane}
\author{Jeremy L. Martin}
\address{Department of Mathematics,
University of California, San Diego,
La Jolla, CA 92093-0112}
\curraddr{Department of Mathematics,
University of Kansas,
Lawrence, KS 66044}
\email{jmartin@math.ku.edu}
\thanks{Supported in part by an NSF Postdoctoral Fellowship.}

\begin{abstract}
Let $m_{12}$, $m_{13}$, \dots, $m_{n-1,n}$ be the slopes of the $\binom{n}{2}$ lines
connecting $n$ points in general position in the plane.  The ideal $I_n$ of all
algebraic relations among the $m_{ij}$ defines a configuration space
called the {\em slope variety of the complete graph}.  We prove that $I_n$
is reduced and Cohen-Macaulay, give an explicit Gr\"obner basis for it, and
compute its Hilbert series combinatorially.  We proceed chiefly
by studying the associated Stanley-Reisner simplicial complex, which has
an intricate recursive structure.
In addition, we are able to answer many questions about the 
geometry of the slope
variety by translating them into purely combinatorial problems concerning
enumeration of trees.
\end{abstract}

\maketitle


\section{Introduction} \label{intro-section}

Let there be given $n\geq 2$ distinct points in the plane, connected in pairs by
$\binom{n}{2}$ lines.  The closure of the locus of all slope vectors
$(m_{12},\dots,m_{n-1,n})$ arising from some such
configuration is an irreducible algebraic variety of dimension $2n-3$,
called the \defterm{affine slope variety of the complete graph} (as we soon explain).
This slope variety turns out to have an unexpectedly rich combinatorial
and geometric structure.  Our techniques for investigating its properties
draw on combinatorics (graph theory and recursive enumeration of trees), commutative
algebra (Gr\"obner bases and Stanley-Reisner theory), and algebraic geometry.
Before stating the main theorem, we set it in context by giving an overview of the theory
of \defterm{graph varieties}, considered by the author in~\cite{JLM1}.

Let $\Pp^2$ be the projective
plane over an algebraically closed field $\fld$,
and let $G$ be a graph with vertices $V$ and edges $E$.  A \defterm{picture} $\Pic$
of $G$ consists of a point $\Pic(v)$ for each vertex, and a line
$\Pic(e)$ for each edge, subject to the conditions that $\Pic(v)\in\Pic(e)$
whenever $v$ is an endpoint of $e$.  Thus the data of $n$ points and $\binom{n}{2}$
lines described earlier is a picture of the complete graph $K_n$ on $n$ vertices.

The set of all pictures of $G$ is called the {\em picture space} $\PS(G)$.
A picture is \defterm{generic} if the points $\Pic(v)$ are all different;
the closure of the locus of generic pictures is called the {\em picture variety}
$\PV(G)$.  This is an irreducible component of $\PS(G)$ of dimension $2|V|$.
Passing to an affine open subset $\PVA(G) \subset \PV(G)$, and projecting onto
an affine space $\Aa_\fld^{|E|}$ whose coordinates correspond to the slopes
of lines $\Pic(e)$,
we obtain the {\em affine slope variety} $\SVA(G)$, of dimension $2|V|-3$.

A \defterm{rigidity circuit} is a graph which admits a decomposition into
two spanning trees, and contains no proper subgraph with that property.  The
most important rigidity circuits are the \defterm{wheels}: a wheel consists
of a cycle with an attached central vertex.
For each rigidity circuit $C$, there is a corresponding \defterm{tree polynomial}
$\TP(C)$, a sum of signed squarefree monomials corresponding to spanning trees appearing
in such decompositions; this polynomial is homogeneous and irreducible.
The affine slope variety
$\SVA(G)$ is cut out set-theoretically in $\Aa^{|E|}$ by the polynomials
$\tau(C)$, where $C$ ranges over all rigidity circuit subgraphs of $G$.
These facts were proven in~\cite{JLM1}.  We can now state the main theorem of this paper.

\vskip 0.1in
\begin{theorem} \label{main-thm} Let $R_n = \fld[m_{12},\dots,m_{n-1,n}]$,
and let $I_n$
be the ideal generated by the tree polynomials of all rigidity circuits in
the complete graph $K_n$.  Then:
\begin{enumerate}
  \item[(i)] The affine slope variety $\SVA(K_n)$ is defined scheme-theoretically by
$I_n$.  That is, $I_n$ is a prime ideal, and $\SVA(K_n) \isom \Spec R_n/I_n$.  
  \item[(ii)] The tree polynomials of the wheel subgraphs of $K_n$
generate $I_n$, and form a Gr\"obner basis with respect to
a certain graded lexicographic order.
  \item[(iii)] $\SVA(K_n)$ has dimension $2n-3$ and degree
    \begin{equation} \label{doublefac}
    \frac{(2n-4)!}{2^{n-2}(n-2)!} = (2n-5)(2n-7) \cdots (3)(1),
    \end{equation}
the number of perfect matchings on $[1,2n-4] = \{1, 2, \dots, 2n-4\}$.
Furthermore, the Hilbert series of $R_n/I_n$ is
    $$\frac{\displaystyle \sum_{k=0}^{n-2} h(n,k) t^k}{(1-t)^{2n-3}},$$
where $h(n,k)$ counts the number of perfect matchings on $[1,2n-4]$ with exactly
$k$ \emph{long pairs}, that is, pairs not of the form $\{i,i+1\}$.
  \item[(iv)] The ring $R_n/I_n$ and the affine slope variety $\SVA(K_n)$ are Cohen-Macaulay.
\end{enumerate}
\end{theorem}
\vskip 0.1in

We begin in Section~\ref{SVA-section} by describing the basic objects---graph varieties
and tree polynomials---in somewhat more detail.  We do this both to make this
paper more self-contained, and because several details of the constructions will
be of importance later on.  The reader is referred to~\cite{JLM1} and~\cite{JLMPHD}
for more leisurely treatments of these subjects.

In the first main part of the paper,
we construct a monomial ideal $J_n$, generated by the initial terms
of tree polynomials $\tau(W)$ of wheels $W \subset K_n$
with respect to a certain graded lexicographic term ordering.
As mentioned previously, the monomials of $\tau(W)$ correspond to \defterm{coupled}
spanning trees of $W$, that is, whose complements are also spanning
trees.  In order to identify the leading term of a wheel polynomial, we need
several technical facts about the valences of vertices in coupled trees;
these facts comprise Section~\ref{valence-section}.  In Section~\ref{leading-tree-section},
we introduce the term ordering and, in Theorem~\ref{initial-ideal-thm},
give a necessary and sufficient combinatorial
condition on monomials (somewhat akin to pattern avoidance in permutations) for
membership in $J_n$.

The second part of the paper consists of Sections~\ref{SR-section}--\ref{hvector-section}.
Here we study the Stanley-Reisner simplicial complex $\SR{n}$
whose faces correspond to squarefree monomials that do not belong to $J_n$.
This simplicial complex has a surprising amount of
combinatorial structure.  First, $\SR{n}$ is pure, and all its
facets may be built up recursively from facets of smaller Stanley-Reisner complexes.
Second, this recurrence may be translated into a bijection between
facets and a combinatorially more natural set, the \defterm{binary total partitions},
which are enumerated by the double factorial numbers~\eqref{doublefac}.
These numbers
therefore give the degree (or multiplicity) of the ideal $J_n$.  Third,
the description of facets leads to a proof that
$\SR{n}$ is shellable and hence Cohen-Macaulay.  This part of the paper
(Section~\ref{shell-section}) is very technical; the
arguments are combinatorially elementary
but do require careful bookkeeping.
The shelling argument leads in turn to a recursive computation of
the $h$-vector of $\SR{n}$; the coefficients of the $h$-vector
enumerate perfect matchings by the number of long pairs (a combinatorial
problem first considered by Kreweras and Poupard~\cite{KP74}).

At this point, we do not yet know that these results on $\SR{n}$ correspond to properties of the 
affine slope variety.  The missing piece is to show that $J_n$ is not too small---precisely,
that it is the initial ideal of
an ideal defining $\SVA(K_n)$ scheme-theoretically.  As it turns out, it is enough to show
that the double factorial numbers give a lower bound for the degree of $\SVA(K_n)$.
We prove this in Sections~\ref{degree-section} and \ref{dpt-section}.
Our approach is to consider a nested family of algebraic subsets of $\SVA(K_n)$ called
 \emph{flattened slope varieties}, whose degree can be bounded from below by a recursive formula
(Theorem~\ref{geom-recurrence-thm}).  We show that this recurrence is equivalent to one
enumerating combinatorial objects called \defterm{decreasing planar trees},
which, like matchings and binary total partitions, are enumerated the double factorials.
Using the fact that $J_n$ is Cohen-Macaulay and has the appropriate codimension and degree,
we conclude that $J_n$ is the (complete) initial ideal of $I_n$ under an appropriate
term ordering; equivalently, the wheel polynomials form a
Gr\"obner basis.  The assertions of the main theorem follow more or less immediately.

Together with \cite{JLM1}, these results constitute the author's doctoral
dissertation \cite{JLMPHD}.  The author thanks his thesis advisor, Mark Haiman,
for his ongoing support.


\section{Preliminaries: Graphs and Tree Polynomials} \label{SVA-section}

We assume that the reader is familiar with the elements of graph theory,
for which a good general reference is~\cite{West}.  We first
fix some notation and terminology.
The symbol $\Nn$ denotes the positive integers.  We abbreviate the set
$\{m, m+1, \dots, n\}$ by $[m,n]$.

A \defterm{graph} $G$ is a pair $(V,E)$, where $V=V(G)$ is a finite set of \defterm{vertices}
and $E=E(G)$ is a set of \defterm{edges}, or unordered pairs of distinct vertices $e=\{v,w\}$.
(Thus we do not allow loops or multiple edges.)  For ease of use, we frequently abbreviate
$\{v,w\}$ by $vw$.
The vertices $v,w$ are the \defterm{endpoints} of $e$.  A \defterm{subgraph} of $G$ is
a graph $G'=(V',E')$ with $V' \subset V$ and $E' \subset E$.  We use the symbols $+$ and
$-$ to denote addition and deletion of edges.

The \defterm{valence} of
$v$ with respect to an edge set $E$, written $\val_E(v)$, is the number of edges in $E$
incident to $v$.  (This is more usually called the \emph{degree}, but we wish to reserve
that term for a different usage.)  A vertex of valence~$1$ is called a \defterm{leaf}.
The \defterm{support} of an edge set $E$ is $V(E) = \{v \st \val_E(v)
> 0\}$.  The \defterm{complete graph} $K_V$ is the graph with vertex set $V$ and
every two vertices adjacent; thus $|E(K_V)|=\binom{|V|}{2}$.
We abbreviate $K_{[1,n]}$ by $K_n$.
For convenience, we frequently ignore the technical distinction between an edge set $E$
and the graph $(V(E),E)$.

Let $v_1,\dots,v_k$ be distinct vertices.  The edge set
$\{v_1v_2,v_2v_3,\dots,v_{k-1}v_k\}$
is called a \defterm{path} from $v_1$ to $v_k$, and if $k\geq 3$ then the edge set
$\{v_1v_2,v_2v_3,\dots,v_{k-1}v_k,v_kv_1\}$
is called a \defterm{cycle} or \defterm{$k$-cycle}.
It is frequently convenient to describe a path or cycle by listing
its vertices in order.

A graph $G$ is \defterm{connected} if every pair of vertices belongs to some common path;
it is a \defterm{tree} if it is connected and contains no cycle.  Equivalently,
a tree may be defined as a connected graph with $|E(G)|=|V(G)|-1$,
or as a graph in which every pair
of vertices belongs to exactly one common path.
A \defterm{spanning tree} of $G$ is a tree $T$
with $V(T)=V(G)$ and $E(T) \subset E(G)$.
The \defterm{connected components} of a graph are its maximal
connected subgraphs.

A graph $G=(V,E)$ is a \defterm{rigidity pseudocircuit} if $E=T \dju T'$,
where $T,T'$ are spanning trees and the symbol $\dju$ denotes
a disjoint union.  A rigidity pseudocircuit is a
\defterm{rigidity circuit} if it contains no other rigidity pseudocircuit as a proper
subgraph.  (For the reasons behind this terminology, see~\cite{GSS}.)
A spanning tree $T\subset E$ is called \defterm{coupled}
if its complement is also a spanning tree;
we denote the set of coupled trees of $G$ by $\SST(G)$.

Let $v_0,v_1,\dots,v_k$ be distinct vertices.
The \defterm{$k$-wheel} $W=W(v_0;\,v_1,\dots,v_k)$ is defined as the graph with edges
$\{v_0v_1,v_1v_2,\dots,v_0v_k\} \cup \{v_1v_2,\dots,v_{k-1}v_k,v_kv_1\}$.
It is easy to check that
every wheel is a rigidity circuit~\cite[Exercise~4.13]{GSS}; the figure below
shows a 2-tree decomposition of $W(v_0;\,v_1,\dots,v_6)$.

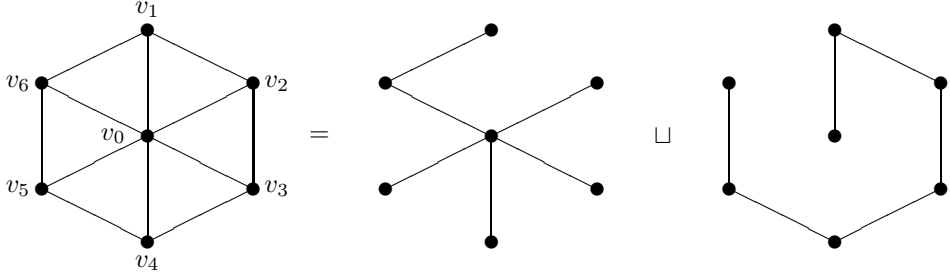
\begin{center}
\begin{picture}(350,110)
\putvertex{ 50, 60} \puttext{37, 60}{$v_0$}
\putvertex{ 50,100} \puttext{50,108}{$v_1$}
\putvertex{ 90, 80} \puttext{99, 80}{$v_2$}
\putvertex{ 90, 40} \puttext{99, 40}{$v_3$}
\putvertex{ 50, 20} \puttext{50, 12}{$v_4$}
\putvertex{ 10, 40} \puttext{ 1, 40}{$v_5$}
\putvertex{ 10, 80} \puttext{ 1, 80}{$v_6$}
\put( 50, 20){\line( 0, 1){80}} 
\put( 10, 40){\line( 2, 1){80}} 
\put( 10, 80){\line( 2,-1){80}} 
\put( 50,100){\line( 2,-1){40}} 
\put( 90, 80){\line( 0,-1){40}} 
\put( 90, 40){\line(-2,-1){40}} 
\put( 50, 20){\line(-2, 1){40}} 
\put( 10, 40){\line( 0, 1){40}} 
\put( 10, 80){\line( 2, 1){40}} 
    \puttext{115,60}{=}
\putvertex{180, 60}
\putvertex{180,100}
\putvertex{220, 80}
\putvertex{220, 40}
\putvertex{180, 20}
\putvertex{140, 40}
\putvertex{140, 80}
\put(180, 20){\line( 0, 1){40}} 
\put(140, 40){\line( 2, 1){80}} 
\put(140, 80){\line( 2,-1){80}} 
\put(140, 80){\line( 2, 1){40}} 
    \puttext{245,60}{$\dju$}
\putvertex{310, 60}
\putvertex{310,100}
\putvertex{350, 80}
\putvertex{350, 40}
\putvertex{310, 20}
\putvertex{270, 40}
\putvertex{270, 80}
\put(310, 60){\line( 0, 1){40}} 
\put(310,100){\line( 2,-1){40}} 
\put(350, 80){\line( 0,-1){40}} 
\put(350, 40){\line(-2,-1){40}} 
\put(310, 20){\line(-2, 1){40}} 
\put(270, 40){\line( 0, 1){40}} 
\end{picture}
\end{center}

The vertex $v_0$ is called the \defterm{center} of $W$, and the vertices $v_1, \dots, v_k$ are
its \defterm{spokes}.  An edge joining two spokes is a \defterm{chord} of the wheel; an edge
joining a spoke to the center is a \defterm{radius}.  We denote the sets of chords and radii
by $\Ch(W)$ and $\Rd(W)$ respectively.  It is notationally convenient to set $v_{k+1}=v_1$,
so that $\Ch(W) = \{ v_iv_{i+1} \st i \in [1,k] \}$.

We now describe certain polynomials associated to rigidity circuits.
Let $n \geq 2$ be an integer and $\fld$ an algebraically closed field.  We will work
over the polynomial ring in $\binom{n}{2}$ variables
    \begin{equation*} \label{define-r-n}
    R_n = \fld[m_{12},\:\dots,\:m_{n-1,n}].
    \end{equation*}
To each edge set $E \subset E(K_n)$ we associate the squarefree monomial
    \begin{equation} \label{mT}
    m_E = \prod_{e \in E} m_e.
    \end{equation}
For every rigidity circuit $C$ in $K_n$, there is an irreducible polynomial, the
\defterm{tree polynomial} of $C$, defined up to sign and having the form 
    \begin{equation} \label{deftreepoly}
    \TP(C) = \sum_{T \in \SST(C)} \varepsilon(T) m_T,
    \end{equation}
where $\varepsilon(T) \in \{1,-1\}$.  By~\cite[Theorem~5.4]{JLM1},
the tree polynomial is irreducible
and homogeneous of degree $|V(C)|-1 = |E(C)|/2$.

We summarize the construction of $\TP(C)$; for more details and examples, see~\cite{JLM1}.
For $1 \leq i
< j \leq n$, regard the edge $\{i,j\}$ as an oriented edge $(i,j)$, and formally set $(j,i)=
-(i,j)$. Fix
a spanning tree $T \subset E(C)$ (not necessarily coupled), and let $\cmp{T} = E(C) \sm T$. For
every edge $e = vw \in \cmp{T}$, the edge set $T+e$ contains a unique cycle
$P_T(e)$, which we may regard as a set of oriented edges $\{(v,w),\,(w,w_1),\, \dots,\,
(w_r,v)\}$.
Let $M$ be the $(|V|-1) \x (|V|-1)$ square matrix with rows indexed by edges $e \in
\cmp{T}$ and columns indexed by the edges $f \in T$, and whose $(e,f)$ entry is
    \begin{equation} \label{deftreepolymatrix}
    M_{e,f} = \begin{cases}
        m_e-m_f & ~\text{if } \phantom{-}f \in P_T(e), \\
        m_f-m_e & ~\text{if } -f \in P_T(e), \\
        0 & ~\text{otherwise.}
    \end{cases}
    \end{equation}
The tree polynomial is then defined as
$\TP(C) = \det M$.
Up to sign, this construction is independent of the choice of the tree $T$.

The \defterm{affine slope variety} $\SVA(n)=\SVA(K_n)$ is defined as follows.
Let $p_1,\dots,p_n$ be $n$ distinct points in the affine plane $\Aa^2_\fld$,
with no two points lying on the same vertical line.  Let $m_{ij} \in \fld$ be the
slope of the unique line joining $p_i$ and $p_j$.  Thus
$(m_{12},\dots,m_{n-1,n})$ is a point in
affine $\binom{n}{2}$-space $\Aa^2_\fld = \Spec R_n$, and $\SVA(n)$
is defined as the closure of the locus of all such points arising from the data
$(p_1,\dots,p_n)$.

By~\cite[Theorem~5.6]{JLM1}, the ideal generated by the tree polynomials
    \begin{equation} \label{define-i-n}
    I_n = \big\langle \TP(C) \st C \subset K_n ~\text{is a rigidity circuit} \big\rangle
    \end{equation}
cuts out the affine slope variety $\SVA(n)$ inside $\Spec R_n$.  In
fact, $I_n$ is generated by the tree polynomials of wheel subgraphs of $K_n$.
We omit the proof of this fact, since we shall 
eventually prove the following stronger result: the wheel polynomials form a Gr\"obner basis for 
$I_n$ with respect to any of a large class of term orderings.


\section{Vertex Valences in Coupled Trees} \label{valence-section}

This section contains several technical facts concerning the valences of
vertices in couples spanning trees of a wheel.  These observations lead eventually
to an explicit identification of the initial monomial of a wheel polynomial $\TP(W)$
with respect to a certain term ordering.

Throughout this section, we fix a $k$-wheel
$W = W(v_0;\,v_1, \dots, v_k) \subset E(K_n)$.
If $T \subset E(W)$ is a coupled spanning tree of $W$, we set $\cmp{T} :=
E(W)\sm T$.

For each $i\in [1,k]$ and each
coupled tree $T \in \SST(W)$, either $\val_T(v_i)=1$ or $\val_T(v_i)=2$.
In addition, not all spokes $v_i$ have the same valence, since $T$ is neither $\Ch(W)$ nor 
$\Rd(W)$.  Thus $\val_T$ may be regarded as a nonconstant function from $[1,k]$ to $[1,2]$.

\begin{lemma} \label{orient-lemma}
Let $T \in \SST(W)$ and $i,j \in [1,k]$.  Then $\cmp{T}$ contains at least one
of the following four edges: $v_0v_i$, $v_0v_j$, $v_iv_{i+1}$, and $v_{j-1}v_j$.
\end{lemma}

\begin{proof}
Suppose not.  Let $i,j$ be a counterexample such that $j-i$ is as small as possible. If
necessary, we may reindex the spokes so that $i \leq j$.  If $j=i$, then $\val_T(v_i)=3$,
which is impossible.  If $j-i=1$, then $T$ contains the cycle $v_0,v_i,v_j,v_0$,
and if $j-i=2$, then $T$ contains the cycle $v_0,v_i,v_{i+1},v_j,v_0$.

Now suppose $j-i > 2$.  Since $T$ contains
the path $v_{i+1},v_i,v_0,v_j,v_{j-1}$, it cannot contain the path
$v_{i+1},v_{i+2},\dots,v_{j-1}$. Let $k$ and $\ell$ be the
least and greatest indices, respectively, such that $v_kv_{k+1} \not\in T$,
$v_{\ell-1}v_\ell \not\in T$, and $i < k \leq \ell < j$. Now $v_0v_k \not\in T$,
otherwise $T$ contains the cycle $v_0,v_i,v_{i+1},\dots,v_k,v_0$.
For a similar reason, $v_0v_\ell \not\in T$.  But
then $k,\ell$ is a counterexample to the lemma,
and $|\ell-k| < |j-i|$, which contradicts the choice of
$i$ and $j$.
\end{proof}

\begin{lemma} \label{another-orient-lemma}
For each $T \in \SST(W)$, at least one of the following conditions is true: either
    \begin{equation*}
    \text{for all}~ i \in [1,k], \quad v_0v_i \in T ~\text{ if and only if }~
    v_iv_{i+1} \in \cmp{T},
    \end{equation*}
or
    \begin{equation*}
    \text{for all}~ i \in [1,k], \quad v_0v_i \in T ~\text{ if and only if }~
    v_{i-1}v_i \in \cmp{T}.
    \end{equation*}
\end{lemma}

\begin{proof}
Suppose that both conditions fail.  That is, there exists $i \in [1,k]$
such that either $v_0v_i,v_iv_{i+1} \in T$ or
$v_0v_i,\,v_iv_{i+1} \in \cmp{T}$.  Interchanging $T$ and $\cmp{T}$ if necessary, we may
assume the former.  Moreover, there exists $j \in [1,k]$ such that either
$v_0v_j,\,v_{j-1}v_j \in T$ or $v_0v_j,\,v_{j-1}v_j \in \cmp{T}$.  The former is ruled
out by Lemma~\ref{orient-lemma}, so the latter must hold; in particular $i\neq j$.

Since $v_0v_i,v_iv_{i+1} \in T$, the edge
$v_0v_{i+1}$ must belong to $\cmp{T}$. 
In particular $\cmp{T}$ contains the path
$v_{j-1},v_j,v_0,v_{i+1}$.  Thus $\cmp{T}$ does not contain the path
$v_{i+1}, v_{i+2}, \dots, v_{j-1}$.
Let $h$ be the largest number in $[i+1,j-1]$ such that
$v_{h-1}v_h \in T$.  Since the path $v_h, v_{h+1}, \dots, v_j,
v_0$ is contained in $\cmp{T}$, the edge $v_0v_h$ must belong to
$T$.  It follows that $T$ contains the edges
$v_0v_i$, $v_iv_{i+1}$, $v_{h-1}v_h$, and $v_0v_h$,
contradicting Lemma~\ref{orient-lemma}.
\end{proof}

\begin{proposition} \label{twotrees-prop}
Let $d: [1,k] \to [1,2]$ be a nonconstant function.  Then there exist exactly two coupled
trees of $W$ for which $\val_T = d$.
\end{proposition}

\begin{proof}
To simplify the notation, write $m_{i,j}$ for $m_{v_iv_j}$.  Also, we put $v_{k+1}=v_1$.
Putting $T=\Rd(W)$ in~\eqref{deftreepolymatrix}, the matrix $M$ becomes
    \begin{equation*}
    \begin{bmatrix}
    m_{0,1}-m_{1,2} & m_{1,2}-m_{0,2} & 0 & \dots & 0 \\
    0 & m_{0,2}-m_{2,3} & m_{2,3}-m_{0,3} & \dots & 0 \\
    \vdots & \vdots & \ddots && \vdots \\
    m_{0,k}-m_{k,1} & 0 & \dots & \dots & m_{k,1}-m_{0,1}
    \end{bmatrix}
    \end{equation*}
Taking the determinant, we obtain
    \begin{eqnarray}
    \TP(W) &=& \prod_{i=1}^k (m_{0,i}-m_{i,i+1}) \;+\;
        (-1)^{k-1} \prod_{i=1}^k (m_{i,i+1}-m_{0,i+1}) \notag \\
    &=& \prod_{i=1}^k (m_{0,i}-m_{i,i+1}) \;-\;
        \prod_{i=1}^k (m_{0,i+1}-m_{i,i+1}). \label{detexp}
    \end{eqnarray}
The monomial $m_{\Ch(W)}$ appears in both products in~\eqref{detexp}, once with coefficient
$+1$ and once with $-1$.  The same is true for the monomial $m_{\Rd(W)}$.  One may
easily verify that
no other cancellation occurs. Accordingly, to enumerate the number of coupled trees by the
valences of spokes, we may substitute $z_i z_{i+1}$ for $m_{i,i+1}$ and $z_i$ for $m_{0,i}$ in
\eqref{detexp} (where the $z_i$ are indeterminates) and change all the $-$'s to $+$'s. This
yields the expression
    \begin{equation*}
    \prod_{i=1}^k (z_i + z_i z_{i+1}) ~+~ 
    \prod_{i=1}^k (z_{i+1} + z_i z_{i+1}) ~-~
    2\left(\prod_{i=1}^k z_i + \prod_{i=1}^k z_i^2 \right) ~=~
    2 \sum_d z_i^{d(i)},
    \end{equation*}
where the sum is taken over all nonconstant functions $d: [1,k] \to [1,2]$.
\end{proof}

\begin{remark}
Proposition~\ref{twotrees-prop} has the following consequence,
which may also be obtained by direct 
counting: every $k$-wheel has exactly $2^{k+1}-4$ coupled spanning trees.
\end{remark}

\begin{defn}
Let $d: [1,k] \to [1,2]$ be a nonconstant function.  The \defterm{type} of a chord
$v_iv_{i+1}$ with respect to $d$ is the pair of numbers $d(i),d(i+1)$.  The \defterm{type} 
of a radius $v_0v_i$ is the pair $d(i-1),d(i+1)$.  If $T \in \SST(W)$, we define 
the type of an edge with respect to $T$ to be its type with respect to $d=\val_T$.
For brevity, we will speak of \emph{type-$11$ chords}, \emph{type-$12$ radii}, etc.
\end{defn}

\begin{lemma} \label{edge-type-lemma}
Let $T \in \SST(W)$, and define the types of chords and radii of $W$ with respect
to the function $\val_T$.  Then:
\begin{enumerate}
  \item[(i)]   Every type-22 chord belongs to $T$.
  \item[(ii)]  Every type-11 chord belongs to $\cmp{T}$.
  \item[(iii)] Every type-22 radius belongs to $\cmp{T}$.
  \item[(iv)]  Every type-11 radius belongs to $T$.
\end{enumerate}
\end{lemma}

\begin{proof}
Let $v_iv_{i+1}$ be a chord.
If (i) fails, then the edges
$v_{i-1}v_i$, $v_{i+1}v_{i+2}$, $v_0v_i$, $v_0v_{i+1}$
all belong to $T$.
If (ii) fails, then those edges all belong to $\cmp{T}$.  In
either case, Lemma~\ref{orient-lemma} is contradicted.

Now let $v_0v_i$ be a radius
Statements (iii) and (iv) are equivalent (switch $T$ and $\cmp{T}$), so we prove only (iii).  
If $\val_T(v_i) = 2$, then $T$ contains the chords $v_iv_{i+1}$ and $v_{i-1}v_i$
by parts (i) and (ii) of the lemma, so the radius $v_0v_i$ belongs to $\cmp{T}$.  On the 
other hand, if
$\val_T(v_i)=1$ and $v_0v_i \in T$, then $\cmp{T}$ contains the chords
$v_iv_{i+1}$ and $v_{i-1}v_i$.  Since $\val_T(v_{i+1}) = \val_T(v_{i-1})=2$ by
hypothesis, the edges
$v_{i-2}v_{i-1}$, $v_{i+1}v_{i+2}$, $v_0v_{i-1}$, $v_0v_{i+1}$
must all belong to $T$, which contradicts Lemma~\ref{orient-lemma}.
\end{proof}

\begin{lemma} \label{yet-another-orient-lemma}
Let $T \in \SST(W)$ and $1 \leq i < j \leq k$.  Suppose that either
    \begin{subequations}
    \begin{equation} \label{suppose1}
    \val_T(v_i) = 1, \qquad \val_T(v_{i+1}) = \dots = \val_T(v_j) = 2, \qquad
    \val_T(v_{j+1}) = 1
    \end{equation}
or
    \begin{equation} \label{suppose2}
    \val_T(v_i) = 2, \qquad \val_T(v_{i+1}) = \dots = \val_T(v_j) = 1, \qquad
    \val_T(v_{j+1}) = 2.
    \end{equation}
    \end{subequations}
Then exactly one of the two chords
$v_iv_{i+1}$, $v_jv_{j+1}$
belongs to $T$.
\end{lemma}

\begin{proof}
Suppose that \eqref{suppose1} holds.  Then the chords $v_{i+1}v_{i+2}$, $v_{i+2}v_{i+3}$,
$\dots$, $v_{j-1}v_j$ all belong to $T$. If both $v_iv_{i+1}$ and 
$v_jv_{j+1}$
belong to $T$, then the path $v_i, v_{i+1}, \dots,
v_{j+1}$ is a connected component of $T$, which is impossible.  On the other hand, if
neither of those chords belong to $T$,
then $v_0v_{i+1}$ and $v_0v_j$ both belong to $T$.  But
then $T$ contains the cycle $v_0,v_{i+1},v_{i+2},\dots,v_j,v_0$, which is impossible,
If we assume \eqref{suppose2} instead of \eqref{suppose1},
the same argument goes through, switching
2 with 1 and $T$ with $\cmp{T}$.
\end{proof}

An alternate formulation of this lemma is as follows.  Let a nonconstant function $d:
[1,k] \to [1,2]$ be given, and let $T$ be a coupled tree with $\val_T=d$.
Traverse the chords of $W$ in order, coloring the type-$12$
chords (of which there are a positive even number) alternately red and blue.  Then
either the red chords all belong to $T$ and the blue chords all belong to
$\cmp{T}$, or vice versa.  Moreover, choosing the color of a single type-$12$ chord suffices to
determine the rest.  Having made such a choice, a radius $v_0v_i$ belongs to $T$
exactly when
$d(i) - \left\vert T \cap \{ v_{i-1}v_i, \, v_iv_{i+1} \} \right\vert = 1$.

Alternatively, if the function $d=\val_T$ is given, then to determine $T$ uniquely
it suffices to specify whether a single type-$12$ radius
$v_0v_i$ belongs to $T$ or to $\cmp{T}$.  Without loss of
generality $v_{i-1}v_i$ is of type $11$ or $22$, and $v_iv_{i+1}$ is of
type $12$.  The value of $d(i)$ determines whether or not $v_iv_{i+1}$ belongs to $T$, 
so the rest of $T$ is determined uniquely as in the preceding paragraph.  That is:

\begin{proposition} \label{tree-type-prop}
Let $T \in \SST(W)$.  Define the type of each edge in $W$ with respect to $\val_T$. Then
$T$ contains all type-$22$ chords, all type-$11$ radii, half the type-$12$ chords, in 
alternation, and a corresponding half of the type-12 radii.
\end{proposition}

We conclude our technical preliminaries
with two results describing the conditions under which a pair of complementary spanning trees may
swap edges.

\begin{lemma} \label{exchange-lemma}
Let $T \in \SST(W)$.  Suppose that $v_{i-1}v_i$ and $v_iv_{i+1}$ belong to
$\cmp{T}$, so that $v_0v_i \in T$.  Assume without loss of generality that
the path in $\cmp{T}$ from $v_i$ to $v_0$ passes through $v_{i+1}$. Then $W$ admits the 
$2$-tree decompositions
    $$E_1 = T - v_0v_i + v_iv_{i+1}, \qquad  E_2 = \cmp{T} - v_iv_{i+1} + v_0v_i$$
and
    $$F_1 = T - v_0v_{i-1} + v_{i-1}v_i, \qquad F_2 = \cmp{T} - v_{i-1}v_i + v_0v_{i-1}.$$
\end{lemma}

\begin{proof}
Clearly $E_1 \dju E_2 = F_1 \dju F_2 = W$.  The edge set $E_1$ is a
tree because $v_i$ is a leaf of $T$, and $E_2$ is a tree because $v_i$ and $v_0$ are in
different connected components of $\cmp{T}-v_iv_{i+1}$.  Meanwhile, $F_1$ is a
tree because $v_i$ and $v_{i-1}$ are in different connected components of $T-v_0v_{i-1}$,
and $F_2$ is a tree because $v_{i-1}$ and $v_0$ are in different connected
components of $\cmp{T}-v_{i-1}v_i$.
\end{proof}

\begin{lemma} \label{another-exchange-lemma}
Let $T \in \SST(W)$. Suppose that $v_{i-1}v_i \in \cmp{T}$
and that $v_0v_i,\,v_iv_{i+1} \in T$, so that
$v_0v_{i+1} \in \cmp{T}$.  Then:

\noindent {\rm (i)} \quad
If $T$ contains at least one radius other than $v_0v_i$, then
$W$ admits the $2$-tree decomposition
    $$E_1 = T - v_0v_i + v_{i-1}v_i, \qquad E_2 = \cmp{T} - v_{i-1}v_i + v_0v_i.$$
\noindent {\rm (ii)} \quad
If $\cmp{T}$ contains at least one radius other than $v_0v_{i+1}$, then 
$W$ admits the $2$-tree decomposition
    $$F_1 = T - v_iv_{i+1} + v_0v_{i+1}, \qquad F_2 = \cmp{T} - v_0v_{i+1} + v_iv_{i+1}.$$
\end{lemma}

\begin{proof}
Clearly $W$ is the disjoint union of $E_1$ and $E_2$ (resp.\ $F_1$ and $F_2$),
so it suffices to show that these edge sets are in fact trees.

\prfpart{i} $E_2$ is a tree because $v_i$ is a leaf of $\cmp{T}$.  If the path
from $v_{i-1}$ to $v_i$ in $T$ does not go through $v_0$, then it must be
$\Ch(W)-v_{i-1}v_i$.  But then $T$ contains at least $k-1$ chords and two radii, which
is impossible.  Therefore $v_{i-1}$ and $v_i$ lie in different connected components of
$T-v_0v_i$, and $E_1$ is a tree.

\prfpart{ii} The path from $v_0$ to $v_{i+1}$ in $T$ is just
$v_0,v_i,v_{i+1}$, so $v_0$ and $v_{i+1}$ belong to different connected components of
$T-v_iv_{i+1}$.  Hence $F_1$ is a tree.  If the path from $v_i$ to $v_{i+1}$ in
$\cmp{T}$ does not go through $v_0$, then it must be $\Ch(W)-v_iv_{i+1}$, which
is impossible.  So $v_i$ and $v_{i+1}$ are in different connected components of
$\cmp{T}-v_0v_{i+1}$.  Hence $F_2$ is a tree.
\end{proof}


\section{The Leading Tree of a Wheel} \label{leading-tree-section}

The main result of this section is Theorem~\ref{initial-ideal-thm}, in which we describe
explicitly the ideal generated by the initial terms of wheel polynomials. Fix once and for all
the following \emph{lexicographic} order $\glex$ on the variables $m_{ij}$:
    \begin{equation*}
    m_{12} \glex m_{13} \glex \dots \glex m_{1n} \glex m_{23} \glex \dots
    \end{equation*}
The corresponding total order for edges of $K_n$ is
    \begin{equation} \label{edgeorder}
    12 \glex 13 \glex \dots \glex 1n \glex 23 \glex \dots
    \end{equation}
We next extend $\glex$ to a term ordering on $R_n$, {\bf graded lexicographic order}, which is
defined as follows: $\prod_{i,j} m_{ij}^{a_{ij}} \glex \prod_{i,j} m_{ij}^{b_{ij}}$ if either
    \begin{equation} \label{glex1}
    \begin{aligned}
    \sum a_{ij} &\:>\: \sum b_{ij}, \qquad \text{or} \\
    \sum a_{ij} &\:=\: \sum b_{ij} \quad \text{and} \quad a_{k\ell} \:>\: b_{k\ell},
    \end{aligned}
    \end{equation}
where $m_{k\ell}$ is the greatest variable (in lexicographic order)
such that $a_{k\ell} \neq b_{k\ell}$.

Associating edge sets with square-free monomials as in \eqref{mT}, we may
regard the term ordering on $R_n$ as defining an extension of the ordering on edges
\eqref{edgeorder} to a total order on subsets of $E(K_n)$.  Then \eqref{glex1} becomes
the following: for $E,F \subset E(K_n)$, $E \glex F$ if either
$|E|>|F|$, or else $|E|=|F|$ and $\max(E \sd F) \in E$, where the symbol $\sd$
denotes the symmetric difference operator.

Given a wheel $W \subset E(K_n)$, we wish to identify the \defterm{leading tree} $LT(W)$
of $W$, that is, the coupled tree of $W$ corresponding to the leading monomial
of $\TP(W)$ (with respect to graded lex order).  We begin by computing the valence of
each spoke of $LT(W)$, using the tools developed in the previous section.  By
Proposition~\ref{tree-type-prop}, this will rule out all but two possibilities for the leading
tree.

\begin{proposition} \label{vertex-valence-prop}
Let $W = W(v_0; v_1, \dots, v_k)$, with $V = V(W) \subset [1,n]$.  Then:

\begin{enumerate}
\item[(i)]  Suppose that $v_0 = \minv(V)$.  Then $\val_{LT(W)}(v_0) = k-1$. 

\item[(ii)] Suppose that $v_0 = \maxv(V)$.  Then $\val_{LT(W)}(v_0) = 1$.

\item[(iii)] Suppose that $v_0 \not\in \{\min(V),\max(V)\}$.  Then for all $i \in [1,k]$,
    \begin{equation*} \label{vdeg-eqn}
    \val_{LT(W)}(v_i) = \begin{cases}
        1 &~ \text{if } v_i \gtv v_0, \\
        2 &~ \text{if } v_i \ltv v_0.
    \end{cases}
    \end{equation*}
\end{enumerate}
In particular, if $v_0$ is the $r$th largest member of $V$,
where $r\in[2,k]$, then
    $$\val_{LT(W)}(v_0) = r-1.$$
\end{proposition}

\begin{proof}
\prfpart{i,ii} Suppose $v_0 = \minv(V)$.  We will show that
if $T$ is a coupled tree with $\val_T(v_0) < k-1$, then $T$ cannot be the leading
tree of $W$.  Note that $\cmp{T}$ contains at least two radii, say $v_0v_i$ and $v_0v_j$.  At least
one of the chords $v_{i-1}v_i$, $v_iv_{i+1}$ belongs to $T$.  If both do, then
by Lemma~\ref{exchange-lemma}, at least one of
    $$T_1 = T - v_{i-1}v_i + v_0v_i, \qquad
      T_2 = T - v_iv_{i+1} + v_0v_i$$
is coupled.  If $v_{i-1}v_i \in T$ and $v_iv_{i+1} \in \cmp{T}$, then
$T_1$ is coupled by Lemma~\ref{another-exchange-lemma}~(i).  But $T_1 \glex T$,
so $T\neq LT(W)$ as desired.  The proof of (ii) is analogous.

\prfpart{iii} Suppose that $v_0$ is neither the minimum nor the maximum element of $V$.
Let $T$ be a coupled tree of
$W$ such that $\val_T(v_i)=1$ for some $v_i \ltv v_0$.  To show that $T \neq LT(W)$, we will
construct a tree $T' \in \SST(W)$ with $T' \glex T$.  There are two cases to consider.

\Case{1}{$v_0v_i \in T$}.  Then $\cmp{T}$ contains the chords $v_{i-1}v_i$ and $v_iv_{i+1}$.
Without loss of generality, we may assume that the path from $v_i$ to $v_0$ in $\cmp{T}$ passes
through $v_{i+1}$.  Then $v_0v_{i-1} \in T$.
By Lemma~\ref{exchange-lemma}, the tree
$T' = T - v_0v_{i-1} + v_{i-1}v_i$ is coupled, and
$v_{i-1}v_i \glex v_0v_{i-1}$, so $T' \glex T$.

\Case{2}{$v_0v_i \in \cmp{T}$}.  Without loss of generality,
$v_iv_{i+1} \in \cmp{T}$ and $v_{i-1}v_i \in T$, so $v_0v_{i+1} \in T$.
Let $T' = T - v_0v_{i+1} + v_iv_{i+1}$.
Note that $T' \glex T$.  If $T'$ is coupled, then we are done.
Otherwise, Lemma~\ref{another-exchange-lemma}~(i) implies that
$v_0v_{i+1}$ is the unique radius in $T$, that is,
$T$ is the path $v_0,\,v_{i+1},\,v_i,\,\dots,\,v_{i-2},\,v_{i-1}$.
Since $v_0\neq\max(T)$, we may choose $j$ such that $v_j \gtv v_0$;
note that $j \neq i$.
Then the tree $T' = T - v_jv_{j+1} + v_0v_{j+1}$ is coupled,
and $T' \glex T$.
\end{proof}

\begin{proposition} \label{center-minmax-prop}
Let $W \subset E(K_n)$ be a $k$-wheel with vertices $V = V(W)$ and center $v_0$.

\begin{enumerate}
\item[(i)] Suppose $v_0 = \minv(V)$.  Label the spokes so that $W = W(v_0; v_1, \dots,
v_k)$ with $v_1 = \maxv\{v_1, \dots, v_k\}$ and $v_2 \gtv v_k$.  Then
    \begin{equation*} \label{crt-min}
    LT(W) = \Rd(W) - v_0v_1 + v_kv_{k+1}.
    \end{equation*}

\item[(ii)] Suppose $v_0 = \maxv(V)$.  Label the spokes so that $W = W(v_0; v_1, \dots,
v_k)$ with $v_1v_2 = \min(\Ch(W))$ and $v_1 \gtv v_2$.  Then
    \begin{equation*} \label{ctr-max}
    LT(W) = \Ch(W) - v_1v_2 + v_0v_2.
    \end{equation*}
\end{enumerate}
\end{proposition}

\begin{proof}
\prfpart{i} By Proposition~\ref{vertex-valence-prop} (i), $LT(W)$ contains exactly one chord.  
Since $v_0v_i \glex v_jv_{j+1}$ for all $i,j$, the unique radius not in $LT(W)$
must be $\min(\Rd(W)) = v_0v_1$.
This implies the desired result because $v_kv_{k+1} \glex v_1v_2$.

\prfpart{ii} For each $1\in[1,k]$, define
    $$T_i = \Ch(W) + v_0v_i - \min\left(v_{i-1}v_i,\,v_iv_{i+1}\right).$$
By Proposition~\ref{vertex-valence-prop}~(ii), $LT(W)$ contains exactly one radius, so $LT(W) =
T_i$ for some $i$.  Note that $T_i = \Ch(W)+v_0v_i-v_1v_2$
for $i=1,2$; in particular $\max(T_1 \sd T_2) = v_0v_2 \in T_2$.  On the other hand,
for $i>2$, we have
    \begin{eqnarray*}
    \max(T_i \sd T_2) &=& \max\left(v_0v_i,\, v_0v_2,\, v_1v_2,\, 
    \min(v_{i-1}v_i, v_iv_{i+1}) \right) \\
    &=& \min(v_{i-1}v_i,\,v_iv_{i+1}) \in T_2.
    \end{eqnarray*}
We conclude that $LT(W)=T_2$.
\end{proof}

In the case that $v_0 \not\in \{\min(V),\max(V)\}$, we have by 
Proposition~\ref{vertex-valence-prop}~(iii)
    \begin{equation*}
    \val_{LT(W)}(i) = \begin{cases}
        1 &~ \text{if } v_i \gtv v_0 \\
        2 &~ \text{if } v_i \ltv v_0.
    \end{cases}
    \end{equation*}

By Proposition~\ref{twotrees-prop}, there are exactly two coupled trees $T,T' \in \SST(W)$
satisfying these conditions.  Moreover
    \begin{eqnarray*}
    T \cap T' &=& \{ \text{chords of type}~22 \} \cup \{ \text{radii of type}~11 \}, \\
    T \sd  T' &=& \{ \text{chords of type}~12 \} \cup \{ \text{radii of type}~12 \}.
    \end{eqnarray*}
Define the \defterm{critical edge} of $W$ to be the maximum element of $T \sd 
T'$.  Thus $LT(W)$ is whichever of $T,T'$ contains the critical edge.

\begin{theorem} \label{initial-ideal-thm}
Let $T$ be a tree with $V(T) \subset [1,n]$.  Then the following are equivalent:
\begin{enumerate}
\item[(i)] There exists a wheel $W \subset K_n$ such that $T = LT(W)$.
\item[(ii)] $T$ contains a path $(v_1, \dots, v_k)$ satisfying the conditions
    \begin{equation} \label{blue}
    \begin{array}{l}
    k \geq 4, \\
    \maxv(v_1, \dots, v_k) = v_1, \\
    \maxv(v_2, \dots, v_k) = v_k, \\
    v_2 \gtv v_{k-1}.
    \end{array}
    \end{equation}
\end{enumerate}
\end{theorem}

\begin{proof}
The easier direction is (ii) $\implies$ (i). Suppose that the path $P = (v_1, \dots, v_k)$
satisfies \eqref{blue}.  Consider the wheel $W = W(v_k; v_1, v_2, \dots, v_{k-1})$.  By
Proposition~\ref{vertex-valence-prop}, $LT(W)$ is a path from $v_k$ to $v_1$.  The two
possibilities for $LT(W)$ are $P$ and the path
    $$P' = (v_k,\,v_2,\,v_3,\,\dots,\,v_{k-2},\,v_{k-1},\,v_1).$$
Since $v_1 \gtv v_k \gtv v_2 \gtv v_{k-1}$, the critical edge of $W$ is
    \begin{eqnarray*}
    \max(P \sd P') &=& \maxv\left(v_1v_2,\, v_{k-1}v_k,\, v_kv_2,\, v_{k-1}v_1\right) \\
    &=& v_{k-1}v_k \in P
    \end{eqnarray*}
so $P = LT(W)$, as desired.

We now show that (i) $\implies$ (ii).  Let $W = W(v_0; v_1, \dots, v_k)$ with $LT(W)=T$.  We
will show that $T$ contains a path $P$ satisfying \eqref{blue}.

\Case{1}{$v_0 = \minv(V)$}.  Reindex the spokes of $W$ so that $v_1 = \maxv(v_1, \dots,
v_k)$ and $v_2 \gtv v_k$.  By Proposition~\ref{center-minmax-prop}~(i) we have
    $$LT(W) = \Rd(W) - v_0v_1 + v_kv_{k+1}.$$
Since $v_0 \ltv v_k \ltv v_2 \ltv v_1$, we may take $P$ to be the path $(v_1,\,v_k,\,v_0,~
v_2)$.

\Case{2}{$v_0 = \maxv(V)$}.  Reindex the spokes of $W$ so that $v_1v_2 =
\min(\Ch(W))$ and $v_1 \gtv v_2$.  By Proposition~\ref{center-minmax-prop}~(ii) we have
    $$LT(W) = \Ch(W) - v_1v_2 + v_0v_2.$$
Let $v_j = \maxv(v_1, \dots, v_k)$.  Obviously $j \neq 2$.  Also $j \neq 3$ (since
$v_1v_2 \llex v_2v_3$, implying $v_1 \gtv v_3$).  Consider the path $P = (v_0,~
v_2,\,v_3,\,\dots,\,v_j)$. The two largest{} vertices of $P$ are $v_0$ and $v_j$.  If $v_2
\ltv v_{j-1}$ then $v_{j-1}v_j \llex v_1v_2$, a contradiction.  So $P$ satisfies
\eqref{blue}.

\Case{3}{$v_0 \not\in \{\min(V),\max(V)\}$}.  Let $e$ be the critical edge of $W$.  By 
definition, $e$ is of type $12$.

\Case{3a}{The critical edge is a chord $v_1v_2$}.  Without loss of generality we may
assume $\val_T(v_1) = 1$ and $\val_T(v_2) = 2$.  Hence $v_1 \gtv v_0 \gtv v_2$.  Now
$v_0v_2 \glex v_1v_2$, so $v_0v_2$ is not of type $12$.  Hence
$v_0v_2$ is of type $11$ and $v_0v_2 \in T$.  Then $v_2v_3 \in \cmp{T}$,
and $d(3)=1$ (otherwise $v_2v_3$ is of type $22$ and $v_2v_3 \in T$, which is
not the case).  Then $v_2v_3$ is of type $12$.  Since $v_1v_2$ is the critical
edge, we have $v_1v_2 \glex v_2v_3$ and $v_1 \ltv v_3$.  In particular $v_1 \neq
\maxv(V)$.  Let $v_j = \maxv(V)$, and let $P$ be the path from $v_1$ to $v_j$ in $T$.  If
$v_0v_j \in T$, then $P$ is the path $(v_1, v_2, v_0, v_j)$, which satisfies
\eqref{blue}.  Otherwise, $P$ is the path
    $$(v_1,\,v_2,\,v_0,\,v_i,\,v_{i-1},\,\dots,\,v_j)$$
for some $i$.  Therefore $\val_T(v_k)=2$ for $i \leq k \leq j-1$, which implies that $v_j$ and 
$v_1$ are respectively the largest{} and second largest{} vertices of
$P$.  Moreover, the chord $v_{j+1}v_j$ is of type $12$, and $v_1v_2$ is the critical 
edge.  So $v_{j+1} \gtv v_2$ and $P$ satisfies \eqref{blue}.

\Case{3b}{The critical edge is a radius $v_0v_1$}.  Without loss of generality we may
assume $\val_T(v_2)=2$ and $\val_T(v_k)=1$.  Hence $v_k \gtv v_0 \gtv v_2$ and $v_0v_1
\llex v_1v_2$.  Thus $v_1v_2$ is not of type $12$ and $\val_T(v_1)=2$.  Let $j$
be the smallest number in $[1,k]$ such that $\val_T(v_j)=1$, and let $P$ be the path $(v_0,
v_1, \dots, v_j)$.  All the edges of $P$ other than $v_0v_1$ and $v_{j-1}v_j$
are chords of type $22$, thus belong to $T$.  So $v_0v_k \in \cmp{T}$ for $2 \leq k
\leq j-1$, and $v_{j-1}v_j \in T$ (because $\val_T(v_{j-1}) = 2$).  Thus $P \subset T$.  
Additionally, $j \geq 3$, and the endpoints of $P$, namely $v_j$ and $v_0$, are respectively
its largest{} and second largest{} vertex.  Finally, $v_{j-1} \gtv v_1$ (because
$v_0v_1$, not $v_0v_{j-1}$, is the critical edge).  Thus $P$ satisfies
\eqref{blue}.
\end{proof}

Clearly, the only path on vertices $[1,4]$ satisfying~\eqref{blue}
is $4213$.  There are two minimal paths on $[1,5]$, namely $53214$ and $52314$.
(Others, such as $53124$, contain subpaths satisfying~\eqref{blue},
in this case $3124$.)  If we let $b(n)$ denote the number of minimal
paths on vertices $[1,n]$ satisfying~\eqref{blue}, then one can check that
the sequence $b(4),b(5),\dots,b(9)$ is $1,2,5,16,61,272$.  This is a subsequence
of the \emph{Euler numbers} or \emph{up/down numbers}, sequence A000111 in~\cite{EIS}.

Many of the results of this section (such as Proposition~\ref{vertex-valence-prop})
hold for all variable orderings which are ``compatible'' with the vertices of $K_n$,
in the sense that if $v,v',v'' \in [1,n]$ and $v'<v''$, then
$m_{vv'} \glex m_{vv''}$.  One might also ask about other natural term orderings on $R_n$,
such as \defterm{reverse lexicographic order}.  This is given by
$\prod m_{ij}^{a_{ij}} \grlex \prod m_{ij}^{b_{ij}}$ if either
$\sum a_{ij} > \sum b_{ij}$, or else
$\sum a_{ij} = \sum b_{ij}$ and $a_{k\ell} < b_{k\ell}$,
where $m_{k\ell}$ is the \emph{smallest} variable for which $a_{k\ell} \neq b_{k\ell}$
(compare~\eqref{glex1}).  Then one can show
that the leading trees of wheels with respect to
reverse lex order are precisely those containing a path $v_1,\dots,v_k$, $k \geq 4$,
such that
$\maxv(v_1, \dots, v_k) = v_1$,
$\maxv(v_2, \dots, v_k) = v_k$, and
$v_2 \ltv v_{k-1}$
(cf.\ Theorem~\ref{initial-ideal-thm}).


\section{The Stanley-Reisner Complex} \label{SR-section}

We begin by recalling the definition of a simplicial complex and some related terminology; for
more detail, see for instance \cite[ch.~5]{BH}. Let $V$ be a finite set of vertices.  An
(abstract)  \defterm{simplicial complex on $V$} is a set $\Delta$ of subsets of $V$ which
contains as members all singleton sets $\{v\}$, $v \in V$, and with the property that if $F
\in \Delta$ and $F' \subset F$, then $F' \in \Delta$.  The elements of $\Delta$ are called
\defterm{faces}.  A maximal face is called a \defterm{facet}.  The \defterm{dimension} of a
face $F$ is $\dim F = |F|-1$, and the dimension of $\Delta$ is $\dim \Delta = \max\{\dim F \st
F \in \Delta\}$.  A simplicial complex is a \defterm{simplex} if it has exactly one facet, and
is \defterm{pure} if all its facets have the same cardinality.

Let $R = \fld[x_1, \dots, x_n]$ and let $J \subset R$ be an ideal generated by squarefree
monomials.  Let $\Delta=\Delta(J)$ be the set of squarefree monomials which do not belong to
$J$. Associating each squarefree monomial $x_{j_1}\dots x_{j_r}$
with the set $\{j_1,\dots,j_r\}$, we may
regard $\Delta$ as a simplicial complex on vertices $[n]$, the \defterm{Stanley-Reisner complex} of
$J$.

The Krull dimension of $R/J$ is $1+\dim\Delta$, and the degree of $R/J$
(that is, its multiplicity as an $R$-module)
equals the number of facets of
$\Delta$.  More particularly, the Hilbert
series of $R/J$ corresponds to the \defterm{$h$-vector} of $\Delta$,
a combinatorial invariant about
which we shall have more to say in Section~\ref{hvector-section}.

We shall study the monomial ideal $J=J_n$ generated by the squarefree
monomials corresponding to leading trees of wheels: that is,
    \begin{equation} \label{define-j-n}
    J_n := \big( m_{LT(W)} \st W \subset K_n ~\text{a wheel} \big) \subset R_n.
    \end{equation}
Thus the Stanley-Reisner simplicial complex $\SR{n}$ of $J_n$ is defined by
    \begin{equation} \label{define-SR-n}
    \SR{n} := \SR{J_n} = \{E \subset E(K_n) \st m_E \not\in J_n\}.
    \end{equation}
Slightly more generally, if $V = \{v_1 < \dots < v_n\}$ is a finite totally ordered set
(typically, $V \subset \Nn$), we define a simplicial complex $\SR{V}$ on $E(K_V)$, isomorphic
to $\SR{n}$, by replacing the edge $ij$ with $v_iv_j$.

Note that the vertices of $\SR{n}$ (in the sense of the definition of a simplicial complex)  
are the edges of $K_n$.  By Theorem~\ref{initial-ideal-thm}, an edge set $E$ is a face of $\SR{n}$
if and only if $E$ contains no path of the form \eqref{blue}.  For this reason, we will say
that a path satisfying \eqref{blue} is \defterm{forbidden}.

\begin{example} \label{SR-example}
No path with fewer than
three edges satisfies \eqref{blue}, so the complexes $\SR{2}$ and $\SR{3}$ are simplices on
$E(K_2)$ and $E(K_3)$ respectively.  There is a unique forbidden path on four vertices,
namely $4213 = \{12,13,24\}$, so the facets of $\SR{4}$ are the edge subsets of $K_4$
that omit one edge from the path 4213, namely, $\{13,14,23,24,34\}$, $\{12,14,23,24,34\}$,
and $\{12,13,14,23,34\}$.  The faces of $\SR{5}$ are those edge sets containing
none of the paths $4213$, $5213$, $5214$, $5314$, $5324$, $53214$, or $52314$.
(These are precisely the paths satisfying~\eqref{blue} which are minimal under
inclusion; see the discussion following the proof of
Theorem~\ref{initial-ideal-thm}.)
\end{example}

The main result on the structure of facets of $\SR{n}$ is as follows.

\begin{theorem} \label{facet-decomp-thm}
Let $n \geq 3$, and let $\SR{n}$ be the Stanley-Reisner simplicial complex
just described.  Then:

\begin{itemize}
\item[(i)] $\SR{n}$ is pure of codimension 
$2n-4$.  Moreover, every facet $F \in \SR{n}$ is $2$-connected,
and may be written uniquely as a disjoint union
    \begin{equation*}
    F ~=~ F^1 \dju F^2 \dju \bs{ 1n }
    \end{equation*}
with the following properties: each $F^i$ is a facet of the complex $\SR{V(F^i)}$, and
    \begin{equation} \label{subfacet-conditions}
    \begin{aligned}
    &V(F^1) \cup V(F^2) = [1,n], \\
    &V(F^1) \cap V(F^2) = \{\max(V(F^1))\}, \\
    &1 \in V(F^1), \\
    &n \in V(F^2).
    \end{aligned}
    \end{equation}

\item[(ii)] Let $F^1,F^2$ be facets of the complexes $\SR{V(F^1)}$ and
$\SR{V(F^2)}$ respectively, satisfying the four conditions just
given.  Then the edge set $F^1 \dju F^2 \dju \{1n\}$ is a facet of $\SR{n}$.
\end{itemize}
\end{theorem}

We will need several facts about connectivity and $2$-connectivity of
graphs, which we summarize here; for details see Chapter~4 of~\cite{West}. Let $G=(V,E)$ be a
connected graph, and denote by $G-v$ the graph obtained by deleting a vertex $v$.
If $G-v$
is disconnected, then $v$ is called a \defterm{cut-vertex} of $G$ (or of $E$).  The vertex $v$
\defterm{separates} the vertices $w$ and $x$ if $w$ and $x$ lie in different connected
components of $G-v$.  Equivalently, $v$ lies on every path between $w$ and $x$.
$G$ is called \defterm{$2$-connected} if it has no cut-vertex.
(In particular, we consider $K_2$ to be $2$-connected, which is not the usual convention.)

A \defterm{block} of $G$ is a maximal $2$-connected subgraph of $G$.  Every edge
of $G$ belongs to exactly one block.  A vertex belongs to more than one block if and only
if it is a cut-vertex; in addition, two blocks share at most one
vertex.  In particular, if $G'=(V',E')$ is a block of $G$, $v \in V'$ is not a cut-vertex of $G$,
and $w \in V \sm V'$, then there is a unique vertex $x \in V'$ such that $x$ is a cut-vertex
of $G$ separating $v$ and $w$.
 
If $v$ is a cut-vertex of a connected graph $G$, then a \defterm{$v$-lobe}
of $G$ is a maximal connected subgraph not having $v$ as a cut-vertex.  Note that
this decomposition is less fine than the block decomposition of $G$: for instance,
if $G$ is a path $v_1,v_2,\dots,v_r$, then each edge is a block, but there
are only two lobes with respect to each cut-vertex.

Finally, we will need the following special case of Menger's Theorem~\cite[p.~167]{West}.
Let $v,w$ be nonadjacent vertices of a graph $G$.  Then exactly one of the
following conditions is true: either $G$ contains a cut-vertex separating $v$ and $w$,
or else $v,w$ lie on a common cycle $C$ in $G$.

We now start our investigation of the simplicial complexes $\SR{n}$.
As noted in Example~\ref{SR-example}, $\SR{2}$ and $\SR{3}$ are simplices;
in particular, $\SR{n}$ is pure of dimension $2n-4$ for
$n \leq 3$.  For $n>3$, the following criterion for nonmembership in $\SR{n}$
will be very useful.

\begin{lemma} \label{nopoly-lemma}
Let $C$ be a cycle on vertices $V \subset [1,n]$.  If $\min(V)$ and $\max(V)$ are not 
adjacent in $C$, then $C$ contains a forbidden path.
\end{lemma}

\begin{proof}
We proceed by induction on $|V|$.  Without loss of generality we may assume $V=[1,n]$.  If
$n=3$, the statement is vacuously true.  If $n=4$, then the only 
possibility is $C = 1,3,4,2,1$,
which contains the forbidden path $4,2,1,3$.

Now suppose that $n>4$, and the vertices $1,n$ are not adjacent.
If $C$ contains the edge $\{1,n-1\}$, then it
contains a forbidden path of the form $n-1,\,1,\,\dots,\,n$.  Otherwise, label
the vertices of $C$ in order as $v_1,\dots,v_n,v_1$,
with $v_1=n$, $v_i=n-1$, and $v_j=1$.  In particular, $1 < i+1 < j < n$.

Let $r = n-i+1$, so that $n > r \geq 4$.  By induction, the $r$-cycle
$v_i,\,v_{i+1},\,\dots,\,v_n,\,v_i$
contains a forbidden path $P$.  If $n-1$ and $v_n$ are not adjacent
in $P$, then $P \subset C$
and we are done.  Otherwise, $n-1$ is an endpoint of $P$, with unique neighbor $v_n$.
Let $P'$ be the path obtained from $P$ by deleting the edge $\{n-1,v_n\}$
and replacing it with $\{n,v_n\}$.  Then $P'$ is also forbidden, and $P' \subset C$ as desired.
\end{proof}

\begin{corollary} \label{minmax-corollary}
Let $F$ be a $2$-connected face of $\SR{n}$, and let $V = V(F)$.  Then
the edge $\{\min(V),\max(V)\}$ belongs to $F$.
\end{corollary}

\begin{proof}
Since $F$ is $2$-connected, it contains a cycle $C$ supported at both
$\min(V)$ and $\max(V)$.  By Lemma~\ref{nopoly-lemma}, these vertices are adjacent in $C$,
hence in $F$.
\end{proof}

By Theorem~\ref{initial-ideal-thm}, no
forbidden path in $K_n$ contains either of the edges
$\{1,n\}$ or $\{n-1,n\}$.  Therefore, both of these edges belong to every facet of
$\SR{n}$.  We will frequently work with the face $\hat{F}=F \sm \{\{1,n\}\}$.
Note that $\hat F$ is a face of $\SR{n}$ of cardinality $|F|-1$.

\begin{lemma} \label{prelim-decomp-lemma}
Let $n \geq 3$, and let $F$ be a facet of $\SR{n}$.  Then both $F$ and $\hat{F}$ are 
connected.  However, $\hat{F}$ is not $2$-connected; in particular, $\hat{F}$ has a 
cut-vertex separating $1$ and $n$.
\end{lemma}

\begin{proof}
Suppose that $F$ is disconnected.  Let $e$ be an edge whose endpoints
are the largest{} vertices of the
connected components of $F$ to which they belong.  Then it
is easy to check that $F+e \in \SR{n}$, which contradicts the hypothesis that 
$F$ is a facet.

Suppose that $\hat{F}$ is disconnected.  Let $H$ be the connected component of $\hat{F}$
containing the edge $\{n-1,n\}$.  Then the vertex $1$ does not belong to $H$;
in particular, $e=\{1,n-1\} \not\in F$.
Let $F' = \hat{F}+e$.  Either $F'$ contains a forbidden path, or it is a face of $\SR{n}$.
If $F'$ contains a forbidden path $P$,
then $P$ contains $e$, hence must be of the form
$n-1,\,1,\,\dots,\,v,\,n$.  But then $P-e \subset F$, which is impossible.  On the other
hand, if $F' \in \SR{n}$, then $F'+\{1,n\} \;=\; F+e \;\in\; \SR{n}$,
which contradicts the hypothesis that $F$ is a facet.  We conclude that $\hat{F}$ is connected.
However, the vertices $1$ and $n$ are not adjacent in $\hat{F}$,
so by Lemma~\ref{nopoly-lemma} they cannot lie in a common cycle.  The last
assertion of this lemma now follows from the special case of Menger's 
Theorem mentioned above.
\end{proof}

\begin{lemma} \label{one-not-cutv-lemma}
Let $n\geq 3$, and let $F$ be a facet of $\SR{n}$.
Then the vertex $1$ is not a cut-vertex of $\hat{F} = F-\{1,n\}$.
\end{lemma}

\begin{proof}
We prove the contrapositive.  Suppose that $1$ is a cut-vertex of $\hat F$.
Let $A$ be the lobe containing
vertex $n$, and let $B = F \sm A$.  Then neither $A$ nor $B$ is empty, and their
only common vertex is $1$.  Let
    \begin{eqnarray*}
    x &=& \maxv(V(B)), \\
    y &=& \minv\{v \in V(A) \st 1v \in \hat{F}\}.
    \end{eqnarray*}
Since $x\not\in V(A)$, we have $x \neq y$.
We consider the cases $x \ltv y$ and $x \gtv y$ separately.
In both cases, the method of proof is to construct a forbidden path in $\hat F$, which is
a contradiction because $\hat{F} \in \SR{n}$.

\Case{1}{$x \ltv y$}.  Suppose that $\hat{F}+xy$ contains a forbidden path
$P$.  Then $xy \in P$, since $\hat F \in \SR{n}$ and $P \not\subset \hat{F}$.  Label
the vertices of $P$ in order as $a, \dots, x,y, \dots, b$.
Then $a,b \geq y > x$, so $a,b \in V(A)$.
Hence the subpath of $P$ from $a$ to $x$ must go through the vertex $1$.  Therefore
$P$ decomposes as
$P_1 \dju P_2 \dju \{xy\} \dju P_3$, where $P_1$ is a path from $a$ to $1$, $P_2$ is 
a path from $1$ to $x$, and $P_3$ is a path from $y$ to $b$.  But then
$P_1 \dju \{1y\} \dju P_3$
is a forbidden path in $\hat{F}$.
It follows that no such path $P$ exists, and $\hat{F}+xy \in \SR{n}$
But then so $F + xy$ is also an element of $\SR{n}$, which
contradicts the hypothesis that $F$ is a \emph{facet} of $\SR{n}$.

\Case{2}{$x \gtv y$}.  The face $\hat{F}$ contains a path
$y,\,1,\,x_1,\,\dots,\,x_r=x$; call this $P$.
If $x_1 \ltv y$, then the subpath of $P$ from $y$ to $x_i$
is forbidden, where $i$ is the smallest index such that $x_i \gtv y$;
this is impossible.  Therefore $x_1 \gtv y$.
On the other hand, $A$ contains a path $P'$ of the form
    $$y=y_1,\,y_2,\,\dots,\,y_s=n$$
with $y_i\neq 1$ for all $i$ (because $1$ is not a cut-vertex of $A$).
Let $t$ be the smallest number such that $y_t>x_1$.
Then $\hat{F}$ contains the forbidden path
    $$x_1,\,1,\,y=y_1,\,y_2,\,\dots,\,y_t,$$
which is a contradiction.
\end{proof}

By this lemma and the remarks following the statement of Theorem~\ref{facet-decomp-thm},
each facet $F$ 
has a unique vertex $a = a(F)$ with the following properties: $a$ is a cut-vertex of $\hat{F}$
separating $1$ and $n$, and $a$ and $1$ belong to the same $2$-connected block of $\hat{F}$.

\begin{lemma} \label{max-cutv-lemma}
Let $F$ be a facet of $\SR{n}$, and let $a=a(F)$ be the vertex just described.
Let $F^1$ be the $a$-lobe of $\hat{F}$ containing $1$, and
let $F^2 = F \sm F^1$ be the union of all other
$a$-lobes.  Then $a = \maxv(V(F^1))$.
\end{lemma}

\begin{proof}
We first prove a weaker statement.  Let $H$ be the unique block of $F$ containing 
$1$ (and thus $a$); we will show that $a = \maxv(V(H))$.  Suppose not, and let
$m = \maxv(V(H)) \gtv a$.  Then $H$ contains the edge $1m$ by 
Corollary~\ref{minmax-corollary}.  Since $H$ is connected, it contains
a path $P_1$ from $1$ to $a$ which does not go through $m$.
Moreover, $F^2$ contains a path from $a$ to $n$, which we can truncate at
the first vertex $v>m$ to obtain a path $P_2$.  Then
$\{1m\} \dju P_1 \dju P_2$
is a forbidden path in $F$ from $m$ to $v$.  This is a contradiction, so $a=m$.

We now prove the full lemma.  Suppose that $m = \maxv(V(F^1)) \gtv a$.  By the weaker
case, $m \not\in V(H)$, so there is a unique cut-vertex $b \in V(H)$ separating $1$
and $m$.  Then $b \neq a$ (otherwise $m=n \in V(F^1)$, which contradicts the definition
of $F^1$) and $b \neq 1$ (by Lemma~\ref{one-not-cutv-lemma}).  Therefore $1 < b < a < m$.
Then $1a \in H$ by Corollary~\ref{minmax-corollary}.  Moreover, $H$
contains a path $P$ from $1$ to $b$ which does not go through $a$.  Additionally,
$F^1 \sm H$ contains a path from $b$ to $m$, which we truncate at the first
vertex greater than $a$ to obtain a path $P'$.
Then $F$ contains a forbidden path, namely $\{1a\} \dju P \dju P'$; this is a
contradiction.
\end{proof}

With all of these technical results in hand, we now proceed to the proof of
the main theorem characterizing facets of the Stanley-Reisner complex $\SR{n}$.

\begin{proof}[Proof of Theorem~\ref{facet-decomp-thm}]
Let $a$, $F^1$ and $F^2$ be as defined in Lemma~\ref{max-cutv-lemma}. To complete the
proof of the theorem, we must show that each $F^i$ is $2$-connected and is a facet of the
complex $\SR{V(F^i)}$, and that the decomposition is unique.  This last assertion is
equivalent to the statement that $a$ is the only cut-vertex of $\hat{F}$.

Suppose that $F^1$ is not a facet of $\SR{V(F^1)}$.  That is, there is some edge $e
\not\in F^1$ such that $F^1 \cup \{e\} \in \SR{V(F^1)}$.  Note that $e \not\in
F$.  Since $F$ is a facet of $\SR{n}$, the edge set $F \cup \{e\}$ must contain a forbidden
path $P$.  Indeed, $P \subset \hat{F} \cup \{e\}$, since no forbidden path contains the edge
$\{1,n\}$. On the other hand, $P \not\subset F^1 \cup \{e\}$ and $P \not\subset F^2$,
since both $F^1 \cup \{e\}$ and $F^2$ are faces of $\SR{n}$. Since both these edge sets
are connected, it follows that $P$ must go through $a$ and have an endpoint $b \in
V(F^1) \sm \{a\}$. But by Lemma~\ref{max-cutv-lemma}, the endpoints of $P$ are not 
its
two greatest vertices. Therefore $P$ is not forbidden, and no such $e$ exists. It follows that
$F^1$ is a facet of $\SR{V(F^1)}$, as desired.  The same argument shows that
$F^2$ must be a facet of $\SR{V(F^2)}$.

We now show by induction on $n$ that $F$ is $2$-connected and has cardinality $2n-3$, and 
that $a$ is the only cut-vertex of $\hat{F}$. The base case $n=3$ is easy: the only facet of
$\SR{3}$ is $F = \{12,13,23\}$.  This has the right cardinality and is $2$-connected, and $F
- 13$ has a unique cut-vertex, namely $2$. Now suppose $n>3$. By induction $F^1$ and
$F^2$ are $2$-connected, and they share the vertex $a$, so $a$ is the only cut-vertex of
$\hat{F}$.  Furthermore, $F^1$ and $F^2$ are precisely the $2$-connected components of
$\hat{F}$.  Since $F = \hat{F}+1n$, and this edge has one endpoint in each of
$V(F^1) \sm \{a\}$ and $V(F^2) \sm \{a\}$, it follows that $F$ is $2$-connected.  
Also, $|V(F^2)| = n-|V(F^1)|-1$, so by induction
    \begin{equation} \label{dim-facet-calc}
    |F| ~=~ 1+\left(2\left|V(F^1)\right|-3\right) +
    \left(2n-2\left|V(F^1)\right|-1\right) ~=~ 2n-3.
    \end{equation}

Finally, suppose that $F^1$ and $F^2$ are facets of $\SR{V(F^1)}$ and
$\SR{V(F^2)}$, respectively, and that the vertex sets $V(F^i)$ satisfy the
conditions \eqref{subfacet-conditions}. So $a = \max(V(F^1))$ is a cut-vertex of
$F^1 \cup F^2$.  Suppose that $F$ contains a forbidden path $P$.  Note that $1n
\not\in P$ and neither $F^1$ nor $F^2$ contains $P$ as a subset, so $P$ must have an
endpoint in $V(F^1) \sm \{a\}$, which contradicts the condition $a =
\max(V(F^1))$. Therefore $F \in \SR{n}$.  Moreover, $|F| = 2n-3$ by
\eqref{dim-facet-calc}, so $F$ is a facet.
\end{proof}

\begin{example} \label{facet-example}
The edge set $F$ shown below is a facet of the complex $\SR{6}$;
this can be checked routinely using the characterization of forbidden paths
in Theorem~\ref{initial-ideal-thm}.  The subsets $F^1$ and $F^2$ are
shown as well.  Note that $4$ is the unique cut-vertex of $\hat{F}=F-16$,
and that the corresponding lobes are exactly $F^1$
and $F^2$.
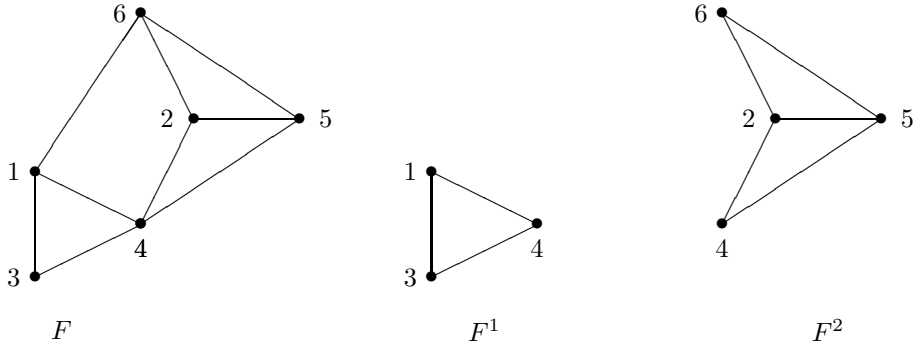
\begin{center}
\begin{picture}(340,130)
  \newsavebox{\faceta}
  \savebox{\faceta}(50,40)[bl]
  {
    \put(10,40){\makebox(0,0){$\bullet$}} \put( 2,40){\makebox(0,0){1}}
    \put(10, 0){\makebox(0,0){$\bullet$}} \put( 2, 0){\makebox(0,0){3}}
    \put(50,20){\makebox(0,0){$\bullet$}} \put(50,11){\makebox(0,0){4}}
    \put(10,40){\line( 0,-1){40}} 
    \put(10,40){\line( 2,-1){40}} 
    \put(10, 0){\line( 2, 1){40}} 
  }
  \newsavebox{\facetb}
  \savebox{\facetb}(80,90)[bl]
  {
    \put(30,50){\makebox(0,0){$\bullet$}} \put(20,50){\makebox(0,0){2}}
    \put(10,10){\makebox(0,0){$\bullet$}} \put(10, 1){\makebox(0,0){4}}
    \put(70,50){\makebox(0,0){$\bullet$}} \put(80,50){\makebox(0,0){5}}
    \put(10,90){\makebox(0,0){$\bullet$}} \put( 2,90){\makebox(0,0){6}}
    \put(30,50){\line(-1,-2){20}} 
    \put(30,50){\line( 2, 0){40}} 
    \put(30,50){\line(-1, 2){20}} 
    \put(10,10){\line( 3, 2){60}} 
    \put(70,50){\line(-3, 2){60}} 
  }
  \puttext{20,0}{$F$}
    \put(0,20){\usebox{\faceta}}
    \put(40,30){\usebox{\facetb}}
    \put(10,60){\line(2,3){40}} 
  \puttext{180,0}{$F^1$}
    \put(150,20){\usebox{\faceta}}
  \puttext{310,0}{$F^2$}
    \put(260,30){\usebox{\facetb}}
\end{picture}
\end{center}

\end{example}


\section{Representing Facets as Binary Trees} \label{rpts-section}

Let $F \in \SR{n}$ be a facet. Iterating the decomposition of Theorem~\ref{facet-decomp-thm},
we can construct a binary tree corresponding to $F$.  In this section, we characterize these
trees explicitly.  Where no confusion can arise, we abbreviate sets of one-digit
integers by a single word, e.g., $13467 = \{1,3,4,6,7\}$.

We begin by reviewing some general facts about planar trees.  A \defterm{rooted tree} is a
tree $T=(V,E)$ with a distinguished vertex $\rt{T}$, the \defterm{root} of $T$.  We
refer to the vertices of a rooted tree as \defterm{nodes};
this is in order to avoid confusion when the nodes are labeled with sets of
vertices of $K_n$.  We
sometimes abuse notation by writing $v \in T$ instead of $v \in V$. If $v,w$
are distinct nodes of $T$,
we say that $v$ is an \defterm{ancestor} of $w$ (equivalently, $w$ is a
\defterm{descendant} of $v$) if $v$ lies on the unique path from $\rt{T}$ to $w$.  If $v$ is
an ancestor of $w$ and $vw \in E$, then $v$ is the \defterm{parent} of $w$ and $w$ is
a \defterm{child} of $v$.  In this case we write $v = v^{(P)}$.  Two vertices with the same
parent are called \defterm{siblings}. We denote by $T|v$ the subtree of $T$ consisting of the
node $v$ and all its descendants.

A \defterm{rooted planar tree} is a rooted tree in which, for every node $v \in V$, the set of
children of $v$ is equipped with a total ordering, the \emph{birth ordering}. That is, we can say
which of two siblings (or two subtrees whose roots are siblings) is \emph{older} than the other.
We denote the $i$th oldest child of $v$ by $v^{(i)}$.  This notation can be iterated: for
instance, $v^{(22)}$ means the second child of the second child of $v$.  A node $w$ is said to
be \defterm{firstborn} if it is the oldest child of its parent; that is,
$w = w^{(P1)}$. If $v = \rt{T}^{(j_1 \dots j_s)}$, then we call the
sequence of numbers $(j_i)$ the \defterm{pedigree} of $v$. A \defterm{binary tree} is a rooted
planar tree in which each node has either zero or two children. The \defterm{traversal}
$\trav(T)$ of a rooted planar tree $T$ is the list of nodes of $T$ in the following order:
first $\rt{T}$, then the nodes of $T^{(1)}$ in traversal order, then the nodes of $T^{(2)}$ in
traversal order, and so on.  (This is the order in which the nodes would be visited in the
course of a depth-first search.)


\begin{defn}
Let $F \in \SR{V}$ be a facet.  Define a binary tree $\TT(F)$ recursively as follows.  If
$|V|=2$, then $\TT(F)$ has one node, labeled by $V$ itself.  Otherwise, $\TT(F)$ is the binary
tree with root $V$, older subtree $\TT(F^1)$, and younger subtree $\TT(F^2)$, where
$F^1$, $F^2$ are as in Theorem~\ref{facet-decomp-thm}.  $\TT(F)$ is called the
\defterm{decomposition tree} of $F$.
\end{defn}

\begin{example} \label{decomp-tree-example}
Let $F$ be the facet shown in Example~\ref{facet-example}. The decomposition tree of $F$ has
root node $123456$, and by our computation of $F^1$ and $F^2$, the children of the root are
$134$ and $2456$.  Continuing in this way, we may calculate the complete decomposition tree,
which is as follows:
    
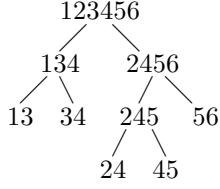
\begin{center}
    \begin{picture}(120,70)
    \put(30,65){\makebox(0,0){123456}}
    \put(15,45){\makebox(0,0){134}}
    \put(50,45){\makebox(0,0){2456}}
    \put( 0,25){\makebox(0,0){13}}
    \put(20,25){\makebox(0,0){34}} 
    \put(45,25){\makebox(0,0){245}}
    \put(70,25){\makebox(0,0){56}}
    \put(35, 5){\makebox(0,0){24}}
    \put(55, 5){\makebox(0,0){45}}
    \put(25,60){\line(-1,-1){10}}
    \put(35,60){\line(1,-1){10}}
    \put(10,40){\line(-1,-1){10}}
    \put(15,40){\line(1,-2){5}}
    \put(50,40){\line(-1,-2){5}}   
    \put(55,40){\line(1,-1){10}}   
    \put(40,20){\line(-1,-2){5}}   
    \put(50,20){\line(1,-2){5}}
    \end{picture}
    \end{center}
\end{example}

There is a bijection $\edg$ from nodes $X$ of $\TT(F)$ to edges of $F$, given by
    \begin{equation*} 
    \edg(X) = \{\min(X),\max(X)\}.
    \end{equation*}
We also write $\edg(T) = \{\edg(X) \st X \in T\}$, where $T$ is a tree or a set of nodes.
The construction of $\TT(F)$ and Theorem~\ref{facet-decomp-thm} immediately
yield the following facts about decomposition trees.

\begin{proposition} \label{admissible-prop}
Let $\TT(F)$ be the decomposition tree of a facet $F \in \SR{V}$.  Then
each node is labeled with at least two vertices of $V$ (for short, $|X|\geq 2$),
and the leaves of $\TT(F)$ are exactly those nodes $X$ for which $|X|=2$.
Each node $X$ that is not a leaf satisfies the following conditions:
    \begin{equation*}
    \begin{array}{ll}
    \text{(a)} & \minv(X) \in X^{(1)} \sm X^{(2)}; \\
    \text{(b)} & \maxv(X) \in X^{(2)} \sm X^{(1)}; \\
    \text{(c)} & X^{(1)} \cup X^{(2)} = X; \\
    \text{(d)} & X^{(1)} \cap X^{(2)} = \{\maxv(X^{(1)})\}.
    \end{array}
    \end{equation*}
\end{proposition}

A tree satisfying these conditions will be called \defterm{admissible}.

\begin{theorem} \label{decomp-tree-thm}
Let $V \subset \Nn$, $|V| \geq 2$.  Let $\Adm(V)$ be the set of all admissible binary trees
with root $V$.
Then the function $\TT$ is a bijection from facets of $\SR{V}$ to $\Adm(V)$.
Moreover, the functions $\TT$ and $\EE$ are inverses.
\end{theorem}

\begin{proof}
Note that $\EE(T) \in \SR{V}$ for all $T \in \Adm(V)$, by Theorem~\ref{facet-decomp-thm}.
Next, we show that the function $\EE$ is injective on $\Adm(V)$.  Let $T \neq T'$ be
admissible trees with the same root $V$.  Then there are nodes $X,X'$ of $T,T'$
respectively which have the same pedigree but different labels.  Choose these two nodes
as close to the root as possible, so that their parents have the same label $Y$.
Then the edge sets $E-\edg(Y)$ and $E'-\edg(Y)$ belong to
different blocks, so the facets $E=\EE(T|Y)$ and
$E'=\EE(T'|Y)$ are distinct and $\EE(T) \neq \EE(T')$ as desired.  It is clear from the
definitions that $\EE(\TT(F))=F$ for every facet $F$, so we are done.
\end{proof}

We may now speak of a \emph{decomposition tree on $V$} as shorthand for a \emph{decomposition 
tree of facets of $\SR{V}$}.

\begin{remark} \label{seniority}
It follows from Proposition~\ref{admissible-prop}
that for each parent node $X$ in a decomposition tree,
one has $X^{(2)} = (X \sm X^{(1)}) \cup \{\maxv(X^{(1)})\}$.  In particular, each older
child determines its younger sibling uniquely, a fact that will be useful later.
\end{remark}


\section{Binary Total Partitions} \label{bpt-section}

By Theorem~\ref{decomp-tree-thm}, the degree of the ideal $J_n$ is the number of decomposition
trees with root $V$.  However, if we try to calculate this number directly using the
conditions of Proposition~\ref{admissible-prop}, we wind up with an awkward recursive formula with
no apparent closed-form solution.  Instead, we construct a bijection from $\Adm(V)$ to a
more easily enumerated set, the \defterm{binary total partitions} of $V$.
We begin by defining these trees
and listing some salient properties.

\begin{defn}
Let $V \subset \Nn$. A \defterm{binary total partition} of $V$ is a binary tree $T$ with
nodes labeled by nonempty subsets of $V$, such that
    \begin{itemize}
    \item $\rt{T}=V$;
    \item If $X$ is a leaf, then $|X|=1$; and
    \item If $X$ has children, then it is their disjoint union, and $\max(X) \in X^{(2)}$.
    \end{itemize}
\end{defn}

It would be simpler to define binary total partitions by ignoring the distinction between the
two children of a given parent.  However, the bijection between decomposition trees and binary
total partitions will be easier to describe if we adopt the convention that each parent and
its younger child have a common maximum.

The set of all binary total partitions of $V$ is denoted $\BT(V)$.  It is elementary to show 
that
    \begin{equation} \label{count-bt}
    |\BT(V)| ~=~ \frac{(2n-4)!}{2^{n-2}(n-2)!} ~=~ (2n-5)(2n-7) \cdots (3)(1).
    \end{equation}
where $n=|V|+1$ \cite[Example~5.2.6]{Stanley-EC2}.

For $T \in \BT(V)$ and $j \not\in V$, define $\Aug_j(T)$ to be the tree obtained by
adding the element $j$ to the root of $T$ (\emph{augmenting $T$ by $j$}).
Clearly $\Aug_j$ is a bijection from $\BT(V)$ to $\BT(V,j) = \{ \Aug_j(T) \st T \in
\BT(V)\}$.  We will say that $T \in \BT(V,j)$ is \defterm{proper} if $j < \min(V)$.

Let $S \subset \Nn$, $j \not\in S$, $j' = \min(S \cup \{j\})$, and $S' = S \cup \{j\} \sm
\{j'\}$.  Define a \emph{straightening map} $\psi_{jj'}$ on trees $T \in \BT(S,j)$ as follows.
For each $k \in S$ with $k<j$, replace all occurrences of $k$ in labels of non-root nodes of
$T$ with the next largest member of $S \cup \{j\}$.  Note that $\psi_{jj'}$ is a bijection
from $\BT(S,j)$ to $\BT(S',j')$; in particular, $\psi_{jj}$ is the identity map. In what
follows, it is frequently convenient and unambiguous to write simply $\psi$ instead of
$\psi_{jj'}$; it is only necessary to specify the subscripts if we want to work with the
inverse function.

\begin{example} \label{phi-example}
Let $S = \{1,2,4,5\}$ and $j = 3$.  Put $j' = \min(S \cup \{j\}) = 1$ and $S' = S \cup \{j\} 
\sm \{j'\} = \{2,3,4,5\}$.  The following figure shows a tree $T \in \BT(S,j)$ and the tree 
$\psi_{jj'}(T) \in \BT(S',j')$.

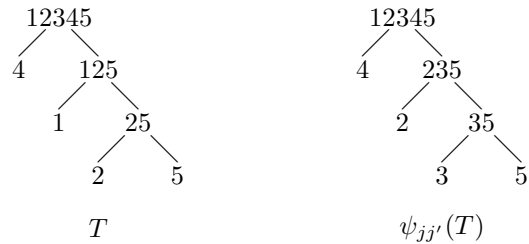
\begin{center}
\begin{picture}(285,100)
\put( 55,90){\makebox(0,0){12345}}
\put( 50,85){\line(-1,-1){10}}
\put( 60,85){\line( 1,-1){10}} 
\put( 40,70){\makebox(0,0){4}}
\put( 70,70){\makebox(0,0){125}}
\put( 65,65){\line(-1,-1){10}}
\put( 75,65){\line( 1,-1){10}}
\put( 55,50){\makebox(0,0){1}}
\put( 85,50){\makebox(0,0){25}}
\put( 80,45){\line(-1,-1){10}}
\put( 90,45){\line( 1,-1){10}}
\put( 70,30){\makebox(0,0){2}}
\put(100,30){\makebox(0,0){5}}
\put( 70,10){\makebox(0,0){$T$}}
\put(185,90){\makebox(0,0){12345}}
\put(180,85){\line(-1,-1){10}}
\put(190,85){\line( 1,-1){10}}
\put(170,70){\makebox(0,0){4}}
\put(200,70){\makebox(0,0){235}}
\put(195,65){\line(-1,-1){10}}
\put(205,65){\line( 1,-1){10}}
\put(185,50){\makebox(0,0){2}}
\put(215,50){\makebox(0,0){35}}
\put(210,45){\line(-1,-1){10}}
\put(220,45){\line( 1,-1){10}}
\put(200,30){\makebox(0,0){3}}
\put(230,30){\makebox(0,0){5}}
\put(200,10){\makebox(0,0){$\psi_{jj'}(T)$}}
\end{picture}
\end{center}
\end{example}

We now construct a bijection between binary total partitions and decomposition trees. Let $S
\subset \Nn$, $j < \min(S)$, and $T \in \BT(S,j)$.  Define a tree $\phi(T)$ recursively as
follows.  If $|S|=1$, then $\phi(T)=T$.  Otherwise, $\phi(T)$ is the tree
    \begin{center}
    \begin{picture}(100,40)
    \put(50,30){\makebox(0,0){$\rt{T}$}}
    \put(45,25){\line(-1,-1){10}}
    \put(55,25){\line( 1,-1){10}}
    \put( 5,10){\makebox(0,0){$\phi(\Aug_j(T^{(1)}))$}}
    \put(95,10){\makebox(0,0){$\phi(\psi(\Aug_m(T^{(2)})))$}}
    \end{picture}
    \end{center}
where $m = \max(\rt{T^{(1)}})$.


\begin{theorem} \label{bin-tot-part-thm}
For every $S \subset \Nn$ and $j < \min(S)$, the function $\phi$ is a bijection 
from $\BT(S,j)$ to $\Adm(S \cup \{j\})$.
\end{theorem}

\begin{proof}
Let $T \in \BT(S,j)$.  Neither augmentation nor applying $\psi$ changes the shape of a
tree, so by induction $\phi(T)$ has the same shape as $T$.  Accordingly, if $X$ is a node
of $T$, we may write $\phi(X)$ for the node of $\phi(T)$ in the same position as $X$.

First, we show that $\phi(T) \in \Adm(S \cup \{j\})$.  If $|S|=1$, this is trivial.
Otherwise, we make the inductive assumption that whenever $|S'|<|S|$ and $j' < \min(S')$
(that is, the members of $\BT(S',j')$ are proper) we have $\phi(T') \in \Adm(S' \cup
\{j'\})$. In particular, the trees $\Aug_j(T^{(1)})$ and $\psi(\Aug_m(T^{(2)}))$ are
proper, so by induction
    $$
    \begin{array}{lllll}
    \phi(T)^{(1)} &=& \phi(\Aug_j(T^{(1)})) &\in& \Adm(X^{(1)} \cup \{j\}) \quad \text{and} \\
    \phi(T)^{(2)} &=& \phi(\psi(\Aug_m(T^{(2)}))) &\in& \Adm(X^{(2)} \cup \{m\}).
    \end{array}
    $$
It is now routine to check that $X = S \cup \{j\}$, $X^{(1)}$ and $X^{(2)}$ satisfy the
conditions of Proposition~\ref{admissible-prop}.

Second, we show that $\phi$ is injective.  If $|S|=1$ or 2 then $|\BT(S,j)|=1$, so there
is nothing to prove.  Otherwise, we assume inductively that $\phi$ maps onto
$\BT(S',j)$ whenever $|S'|<|S|$. Let $T,U \in \BT(S,j)$ with
$\phi(T)=\phi(U)$.  Then
    $$\phi(\Aug_j(T^{(1)})) ~=~ \phi(T)^{(1)} ~=~ \phi(U)^{(1)} ~=~ 
    \phi(\Aug_j(U^{(1)})).$$
By induction, we can cancel the $\phi$'s from the outer terms, and $\Aug_j$ is
a bijection, so $T^{(1)} = U^{(1)}$.  Now, setting $m=\max(\rt{T^{(1)}}) =
\max(\rt{U^{(1)}})$, we have
    $$\phi(\psi(\Aug_m(T^{(1)}))) ~=~ \phi(T)^{(2)} ~=~ \phi(U)^{(2)} ~=~ 
    \phi(\psi(\Aug_m(U^{(1)}))).$$
As before, it follows from the inductive hypothesis and the bijectivity of
$\psi$ and $\Aug_m$ that $T^{(2)}=U^{(2)}$.  Of course $\rt{T} = \rt{U} = S \cup \{j\}$, 
so $T=U$ as desired.

Finally, we show that $\phi$ is surjective.  If $|S|=1$ or $2$, then $|\Adm(S \cup \{j\})|=1$,
so there is nothing to prove.  Otherwise, we make the inductive assumption that $\phi$ is
surjective on $\BT(S',j)$ whenever $|S'|<|S|$.  Suppose that $j < \min(S)$ and $U \in \Adm(S
\cup \{j\})$.  Note that $U^{(1)} \in \Adm(S^{(1)} \cup \{j\})$ and $U^{(2)} \in \Adm(S^{(2)}
\cup \{m\})$, where $S^{(1)}$ and $S^{(2)}$ are nonempty disjoint sets with $S^{(1)} \cup
S^{(2)} = S$, $m = \max(S^{(1)})$, and $\max(S) \in S^{(2)}$.  Let $k = \min(S^{(2)} \cup
\{m\})$ and
    $$
    \begin{array}{lllll}
    T^{(1)} &=& \Aug_j^{-1}(\phi^{-1}(U^{(1)})) &\in& \BT(S^{(1)}), \\
    T^{(2)} &=& \Aug_m^{-1}(\psi_{mk}^{-1}(\phi^{-1}(U^{(2)}))) &\in& \BT(S^{(2)}).
    \end{array}
    $$
Let $T$ be the tree
    \begin{center}
    \begin{picture}(100,30)
    \put(50,20){\makebox(0,0){$S \cup \{j\}$}}
    \put(45,15){\line(-1,-1){10}}
    \put(55,15){\line( 1,-1){10}}
    \put(35,0){\makebox(0,0){$T^{(1)}$}}
    \put(65,0){\makebox(0,0){$T^{(2)}$}}
    \put(100,0){\makebox(0,0){.}}
    \end{picture}
    \end{center}
Then $\phi(T)=U$ and we are done.
\end{proof}

For convenience, we summarize here our present knowledge about the
simplicial complex $\SR{n}$.  Note that $\dim\,I_n = 
\dim\,J_n = 2n-3$ and $\ini(I_n) \supset J_n$, so $\deg\,I_n =
\deg\,\ini(I_n) \leq \deg\,J_n$.  Therefore:

\begin{theorem} \label{facet-count-thm}
The function $\theta = \phi \circ \Aug_1$ is a bijection from
$\BT([2,n])$ to $\Adm([1,n])$.  In particular,
$\SR{n}$ has exactly
$$\displaystyle\frac{(2n-4)!}{2^{n-2}(n-2)!}$$
facets.  Moreover, this number equals the degree of the ideal $J_n$, and is
an upper bound for the degree of $I_n$.
\end{theorem}


\section{A Shelling of $\SR{n}$} \label{shell-section}

We prove in this section that the simplicial complex $\SR{n}$ is shellable, hence
Cohen-Macaulay, for all $n \geq 2$.  We begin by recalling briefly the definition of
shellability.  (There are many equivalent definitions; see, e.g., \cite[pp.~214--219]{BH}.)

\begin{defn}
Let $\Delta$ be a pure simplicial complex. A total ordering $\gfac$ of the facets of $\Delta$
is called a \defterm{shelling} if each facet $F$ has a subset $\shs(F)$ such that

(SH1) If $F \gfac G$, then $\shs(F) \not\subset G$, and

(SH2) For each $e \in \shs(F)$, there is some $G \lfac F$ such that $F \sm G =
\{e\}$.
\end{defn}

We first fix some notation.
Let $T$ be a decomposition tree.  Denote by $\LL(T)$ the set consisting of the root of $T$
together with all firstborn nodes.  For $F$ a facet of $\SR{n}$, we set $\LL(F) =
\LL(\TT(F))$.

Order the finite subsets of $\Nn$ as follows: $X \gint Y$ if $|X| > |Y|$, or if $|X| = |Y|$
and $\min(X \sd Y) \in Y$.  Let $F,G \in \SR{n}$ be facets, with traversals
    $$\begin{array}{ccc}
    \trav(\TT(F)) &=& (X_1, \; \dots, \; X_N), \\
    \trav(\TT(G)) &=& (Y_1, \; \dots, \; Y_N),
    \end{array}$$
where $N = 2n-3$.  Let $k = \min\{i \st X_i \neq Y_i\}$.  Then
we define a total order on facets by putting $F \gfac G$ if $X_k \gint Y_k$.  Note that in 
this case $X_k$ and $Y_k$ have the same pedigree.

We will prove that this ordering satisfies the conditions {\bf (SH1)} and {\bf (SH2)}.
First, we prove a fact to be used in the proof.

\begin{lemma} \label{span-decomp-lemma}
Let $T$ be a decomposition tree on $V \subset \Nn$.  Then the edge set $\edg(\LL(T))$ is a
spanning tree of the complete graph on $V$.
\end{lemma}

\begin{proof}
Suppose that $\edg(\LL(T))$ contains a cycle $C$.
Let $Y$ be the least common ancestor of the nodes in $T$ corresponding
to edges in $C$, and let $e=\edg(Y)$.
Then the endpoints of $e$ are the minimum and maximum vertices of $C$,
so $e$ lies on $C$ by Lemma~\ref{nopoly-lemma}.
Let $e'$ be the other edge of $C$ incident to the vertex
$\maxv(V(C)) = \maxv(Y)$.  The node $X$ corresponding to $e'$ belongs to the \defterm{right
rib} of $T|Y$, that is, the set
    \begin{equation} \label{rib}
    \rr(T) = \left\{ \rt{T}, ~ \rt{T}^{(2)}, ~ \rt{T}^{(22)}, ~ \dots \right\}.
    \end{equation}
In particular $X \not\in \LL(T)$, which is a contradiction.  Therefore $\edg(\LL(T))$ does not
contain any cycles.  Since its cardinality is $n-1$, it is a spanning tree of $V$.
\end{proof}

\begin{theorem} \label{shelling-thm}
The order $\gfac$ is a shelling order on the facets of $\SR{V}$, with $\shs(F)$ defined
recursively by
    \begin{equation*} \label{shellset}
    \shs(F) =
    \begin{cases}
      \emptyset & ~\text{if } |V| \leq 3\\
      \shs(F^2) & ~\text{if } V(F^1) = \{\min(V),\tmin(V)\}\\
      \shs(F^2) \cup \LL(F^1) & ~\text{if } {\rm otherwise\/}
    \end{cases}
    \end{equation*}
where $\tmin(V)$ denotes the second smallest element of $V$.
\end{theorem}

\begin{proof}
We may assume without loss of generality that $V = [1,n]$.  In the first part of the proof, we
will show that the order $\gfac$ satisfies {\bf SH1}.  Let $F,G$ be facets of $\SR{n}$ such
that $F \gfac G$.  Put $T = \TT(F)$, $U = \TT(G)$, $\trav(T) = (X_1, \dots, X_N)$, and
$\trav(U) = (Y_1, \dots, Y_N)$, with $X_i=Y_i$ for $i<k$ and $X_k \gint Y_k$.  
For ease of use, we abbreviate $X=X_k$ and
$Y=Y_k$.

The parent nodes $X^{(P)}=Y^{(P)}$ have the same label, a set of cardinality $\geq 4$.  If $X$ and
$Y$ are the younger siblings of their parents,
then $X^{(P1)}=Y^{(P1)}$ by definition of $k$, but then $X=Y$ by
Remark~\ref{seniority}, a contradiction. Therefore $X=X^{(P1)}$ and $Y=Y^{(P1)}$.
In particular
$E \subset \shs(F)$, where $E = \edg(\LL(T|X))$.  Now $X \gint Y$, so in particular $X \neq
\{\min(X^{(P)}),\tmin(X^{(P)})\}$. Therefore $\edg(X) \in E$.

We claim that $E \not\subset G$.  It is sufficient to prove that $E \not\subset
\edg(U|Y^{(P)})$, because if $Z \in U$ is not a descendant of $Y^{(P)}$, then the endpoints of
$\edg(Z)$ do not both belong to $Y^{(P)}=X^{(P)}$, so $\edg(Z) \not\in E$.

Suppose first that $\max(X) > \max(Y)$.  Then $\max(X) \not\in Z$ for all $Z \in U|Y$;
$\min(X) \not\in Z$ for all $Z \in U|Y^{(P2)}$; and $\max(X) < \max(Y^{(P)})$.  Hence $\edg(X) 
\not\in \edg(U|Y^{(P)})$, establishing the claim.

Now, suppose that $\max(X) \leq \max(Y)$.  Let $v \in X \sm Y$. By
Lemma~\ref{span-decomp-lemma}, $E$ is a spanning tree of the vertices in $X$.  In particular,
$E$ contains a path $P$ from $v$ to $\min(X)=\min(Y)$.  If $\max(X)=\max(Y)$, then $\max(Y)$
is a leaf of $E$ (since $\edg(X)$ is the only edge having it as an endpoint), while if
$\max(X) < \max(Y)$, then $\max(Y) \not\in V(E)$.  In either case $\max(Y) \not\in
V(P)$.  Note that $v \in Y^{(P2)} \sm Y$ and $\min(Y) \in Y \sm Y^{(P2)}$. By
Theorem~\ref{facet-decomp-thm}, $\max(Y)$ is a cut-vertex of $\edg(U|Y^{(P)}) \sm
\{\edg(Y^{(P)})\}$.  So every path from $\min(Y)$ to $v$ in $\edg(U|Y^{(P)})$ must either pass
through $\max(Y)$ or include the edge $\edg(Y^{(P)})$.  But $P$ does neither of these things,
so $P \not\subset \edg(U|Y^{(P)})$, establishing the claim and completing the proof that
$\gfac$ satisfies {\bf SH1}.

\vskip 0.1in

For the second part of the proof, let $F \in \SR{n}$ be a facet, $T = \TT(F)$, $e \in
\shs(F)$, and $X$ the node of $T$ corresponding to $e$.  We will show that $\SR{n}$ contains a
facet $G$ satisfying condition {\bf SH2}.  We proceed by induction on $n$; there is nothing to
prove if $n \leq 3$.

Suppose first that $e \in \shs(F^2)$.  By induction, the simplicial complex
$\SR{V(F^2)}$ has a facet $G^2 \lfac F^2$ such that $F^2 \sm G^2 =
\{e\}$.  Therefore, the facet of $\SR{n}$ given by $
    F^1 \cup G^2 \dju \big\{\{\min(V(F)),\max(V(F))\}\big\}$
satisfies {\bf SH2}.

On the other hand, suppose that $V(F^1) \neq \{1,2\}$ and $e \in \LL(F^1)$.  We 
consider two separate cases: either $X$ has children or it does not.

\Case{1}{$X$ has children}.  Let $Y = X^{(P)}$, $E = \edg(T|Y)$,
$F_1 = \edg(T|X^{(1)})$, and $F_2 = \edg(T|X^{(2)}) \cup \edg(T|Y^{(2)})$.
Note that $F_1$ is a facet of $\SR{X^{(1)}}$, and $F_2 \in \SR{X^{(2)} \cup
Y^{(2)}}$.  Moreover, since $X^{(2)} \cap Y^{(2)} = \{\max(X^{(2)})\}$, we have
    $$|F_2| ~=~ (2|X^{(2)}|-3) + (2|Y^{(2)}|-3) ~=~ \dim(\SR{X^{(2)} \cup Y^{(2)}})$$
so $\SR{X^{(2)} \cup Y^{(2)}}$ has a facet $F'_2$ of the form $F_2\cup\{e'\}$.  
Let $G' = F_1 \cup F'_2 \cup \{\edg(Y)\}$; this is a facet of $\SR{Y}$ because
    $$X^{(1)} \cup (X^{(2)} \cup Y^{(2)}) = Y \quad\text{and}\quad
    X^{(1)} \cap (X^{(2)} \cup Y^{(2)}) = \{\max(X^{(1)})\}.$$
Let $G$ be the facet of $\SR{V}$ obtained by replacing $T|Y$ with $\TT(G')$.  Then
$F \sm G = \{e\}$.  Moreover, the traversals of $\TT(F)$ and $\TT(G)$ first differ at 
$Y^{(1)}$, which is $X$ in $\TT(F)$ and $X^{(1)}$ in $\TT(G)$.  Since $X^{(1)}
\subsetneq X$, we have $F \gfac G$, so $G$ satisfies {\bf SH2}.

\Case{2}{$X$ is a leaf}.  In particular $|X|=2$.  If $X \neq \{\min(X^{(P)}),
\tmin(X^{(P)})\}$, then replacing $X$ by that set produces a facet $G$ which is easily seen to
satisfy {\bf SH2}.  Now suppose that $X = \{\min(X^{(P)}), \tmin(X^{(P)})\}$.  Put $X_0 = X$
and $X_i = X_{i-1}^{(P)}$ for $i \geq 1$.  Notice that $T^{(P)}$ is not in the right rib of
$T$; if it were, it would have to be the root of $T$ (since $X$ is itself in the left subtree
of $T$), which would imply that $X$ is the only node in the left subtree, hence not a member
of $\shs(F)$.  Accordingly, let $Y$ be the youngest ancestor of $X$ such that $X \in
T|Y^{(1)}$; then $Y = X_s$ for some $s \geq 2$.  (If no such $Y$ exists, then $X^{(P)} \in
\rr(T)$, but $X \in \LL(T)$ is in the left subtree of $\rt{T}$, so that $X^{(P)} = \rt{T}$.)
Also, let $Z = X^{(P2)}$.  That is, $T|Y$ has the following form:
    
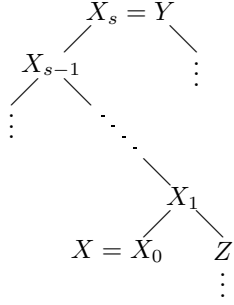
\begin{center}
    \begin{picture}(170,120)
    \put(50,110){\makebox(0,0){$X_s=Y$}}
    \put(35,105){\line(-1,-1){10}}
    \put(65,105){\line(1,-1){10}}
    \put(20,90){\makebox(0,0){$X_{s-1}$}}
    \put(75,90){\makebox(0,0){$\vdots$}}
    \put(15,85){\line(-1,-1){10}}
    \put(25,85){\line( 1,-1){10}}
    \put(5,70){\makebox(0,0){$\vdots$}}
    \put(39,71){\line(1,0){1}}
    \put(43,67){\line(1,0){1}}
    \put(47,63){\line(1,0){1}}
    \put(51,59){\line(1,0){1}}
    \put(55,55){\line( 1,-1){10}}
    \put(70,40){\makebox(0,0){$X_1$}}
    \put(65,35){\line(-1,-1){10}}
    \put(75,35){\line( 1,-1){10}}
    \put(45,20){\makebox(0,0){$X=X_0$}}
    \put(85,20){\makebox(0,0){$Z$}}
    \put(85,10){\makebox(0,0){$\vdots$}}
    \end{picture}
    \end{center}
For $i \in [s-1]$, define
$E_i = \edg(T|X_i) \sm \{e\}$.
Each $\edg(T|X_i)$ is a facet of $\SR{X_i}$, hence $2$-connected by
Theorem~\ref{facet-decomp-thm}, so every $E_i$ is connected. On the other hand, we claim that
$E_i$ is not $2$-connected. More specifically, we claim that for every $i$, the vertex
$m=\max(X_{s-1})$ is a cut-vertex separating $\min(X)$ and $\max(X)$. If $i=1$, the vertex
$\min(X)$ is incident to exactly two edges in $\edg(T|X_1)$, namely $\edg(X)=e$ and $\edg(X_1)
= \{\min(X_1),\max(X_1)\} = \{\min(X),m\}$, so $m$ is the only neighbor of $\min(X)$ in $E_1$.  

Now suppose that $1 < i \leq s-1$.  Note that
both $\edg(X_1) = \{\min(X),m\}$ and
$\edg(Z) = \{\max(X),\max(Z)=m\}$ belong to $E_i$.
Suppose that $E_i$ contains a path $P$ from $\min(X)$ to $\max(X)$ which does not pass through
$m$.  Let $e'$ be the edge incident to $\min(X)$ in $P$, and $W$ the corresponding node of
$T$.  Then $W$ is not one of the nodes $X_j$, since $m =
\max(X_j)$ for all $j$.  Also, $W \not\in T|Z$, because $\min(X) \not\in Z$.  Therefore $W \in
T|X_j^{(1)}$ for some $j \in [2,i]$.  In particular $\min(X) \in X_j^{(1)} \cap
X_{j-1}$, so $\min(X) = \max(X_j^{(1)})$ and $\min(X) \neq \min(V(P))$.  But then $P \cup
\left\{ \edg(X_1), \edg(Z) \right\}$ is a cycle in $E_i$ in which the smallest vertex,
namely $\min(V(P))$, and the largest vertex, namely $m$, are not adjacent.
This contradicts Lemma~\ref{nopoly-lemma} and completes the proof of the claim.

Let $F_1$ be the block of $E_{s-1}$ containing the vertex $\min(Y)$, and let
$F_2 = (E_{s-1} \sm F_1) ~\cup~ \edg(T|Y^{(2)})$.  Also let $V_i = V(F_i)$ for
$i=1,2$.  Note that
    \begin{subequations}
    \begin{eqnarray}
    V_1 \cup V_2  &=& Y, \label{note1} \\
    V_1 \cap V_2  &=& \{m\}, \qquad \text{and} \label{note2} \\
    V_1 &\subsetneq& X_{s-1}. \label{note3}
    \end{eqnarray}
    \end{subequations}
For $i=1$ or $2$, let $F'_i$ be a facet of $\SR{V(F_i)}$ containing $F_i$.  Now
    $$|F_1| + |F_2| = |F_1 \cup F_2| = |\EE(T|Y) \sm \{\edg(X),\edg(Y)\}| = 2|Y|-5$$
and
    $$|F'_1| + |F'_2| = (2|V_1|-3) + (2|V_2|-3) = 2|Y|-4,$$
the last equality following from \eqref{note1} and \eqref{note2}. Therefore
$F'_1 \cup F'_2 = F_1 \cup F_2 \cup \{e'\}$, for some edge $e' \neq e$, and the face
$G' = F'_1 \cup F'_2 \cup \{\edg(Y)\}$
is a facet of $\SR{Y}$ with $\edg(T|Y) \sm G' = \{e\}$. Let $G$ be the facet of
$\SR{V}$ obtained by replacing $T|Y$ with $\TT(G')$. Then $F \sm G = \{e\}$ as well.  
Furthermore, \eqref{note3} implies that $F \gfac G$.  Therefore $G$ is the desired facet
satisfying {\bf SH2}.
\end{proof}


\section{The $h$-vector of $\SR{n}$} \label{hvector-section}

Since $\SR{n}$ is shellable, its $h$-vector
    $$h(\SR{n}) ~=~ (h(n,0), \; h(n,1), \; \dots)$$
has the following combinatorial interpretation~\cite[Corollary~5.1.14]{BH}: for any shelling
order $\gfac$, $h(n,k)$ is the number of facets $F \in \SR{n}$ with $|\shs_{\gfac}(F)|=k$.
In this section, we prove that the numbers $h(n,k)$ have another, more elementary
combinatorial interpretation: they enumerate perfect matchings by the number of
\defterm{long pairs}.  These numbers were first investigated by
Kreweras and Poupard~\cite{KP74}; we discovered the connection using the
On-Line Encyclopedia of Integer Sequences~\cite{EIS}.
We begin by obtaining a recurrence for $h(n,k)$, using the description of facets
and the shelling order of Theorem~\ref{shelling-thm}.  We then show that the
recurrence is equivalent to one enumerating matchings by the number of long pairs

Let $F \in \SR{n}$ be a facet and $T = \TT(F)$ the corresponding decomposition tree.  Define
$\shs(T) ~=~ \left\{ X \in T \st \edg(X) \in \shs(F) \right\}$.

Let $\theta: \BT{[2,n]} \to \Adm([1,n])$ be the bijection of Theorem~\ref{facet-count-thm}.  
Recall that $T$ and $\theta(T)$ have the same shape, so we may define $\shs(U)$ to be the set
of nodes in the same positions as the nodes in $\shs(\theta^{-1}(U))$.

Recall the definition \eqref{rib} of the right rib $\rr(T)$ of a binary tree.
As before, denote the second smallest element of a set $S$ by $\tmin(S)$.

\begin{lemma} \label{hvector-lemma-one}
Let $F \in \SR{n}$ be a facet, $T = \TT(F)$ the corresponding decomposition tree, and $X$ a
firstborn node of $T$.  Then:

\begin{enumerate}
\item[(i)] $X \not\in \shs(T)$ if and only if
$X = \left\{ \min(X^{(P)}), \tmin(X^{(P)}) \right\}$ and
$X^{(P)} \in \rr(T)$.

\item[(ii)] $\theta(X) \not\in \shs(\theta(T))$ if and only
if $\theta(X) = \{ \min(\theta(X)^{(P)}) \}$ and
$\theta(X)^{(P)} \in \rr(\theta(T))$.
\end{enumerate}
\end{lemma}

\begin{proof}
\prfpart{i} The first condition implies that $X \not\in \shs(X^{(P)})$; the second condition
implies that $\shs(X^{(P)}) = \shs(F) \cap \edg(T|X^{(P)})$.  Together, they imply that $X
\not\in \shs(F)$.  On the other hand, suppose that either or both conditions fail. Let $Y$ be
the youngest ancestor of $X$ which does not belong to $\rr(T)$.  Then $Y$ is firstborn, and $X
\in \shs(T|Y^{(P)}) \subset \shs(F)$.

\prfpart{ii} This is obtained by translating the conditions from Lemma~\ref{hvector-lemma-one}(i)
into conditions on $\theta(T)$, using the definition of $\theta$.
\end{proof}

Let $T \in \BT{[2,n-1]}$ and $X \in T$.  Define a tree $\gamma(T,X) \in \BT{[1,n-1]}$ as
follows: replace $T|X$ with the tree
    
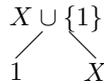
\begin{equation*}
    \begin{picture}(40,30)
    \put(20,25){\makebox(0,0){$X \cup \{1\}$}}
    \put(15,20){\line(-1,-1){10}}
    \put(25,20){\line( 1,-1){10}}
    \put( 5, 5){\makebox(0,0){$1$}}
    \put(35, 5){\makebox(0,0){$X$}}
    \end{picture}
    \end{equation*}
and append $1$ to the label of every ancestor of $X$.  The map
    $$\gamma:~ \left\{ (T,X) \st T \in \BT{[2,n-1]}, ~ X \in T \right\}
    ~\to~ \BT{[1,n-1]}$$
is clearly one-to-one, and these sets have the same cardinality, so $\gamma$ is a bijection.
By incrementing all numbers in the labels of all nodes of $\gamma(T,X)$, we obtain a bijection
$\tilde\gamma$ which maps pairs $(T,X)$ as above to $\BT{[2,n]}$.  
Moreover,
    \begin{equation} \label{effectshell}
    \shs \left( \tilde\gamma(T,X) \right) ~=~ \begin{cases}
        \shs(T) & ~\text{if } X = \rt{T}, \\
        \shs(T) \cup \{1,X\} & ~\text{if } X ~\text{is firstborn and}~ X \not\in \shs(T), \\
        \shs(T) \cup \{1\} & ~\text{otherwise,}
    \end{cases}
    \end{equation}
where $\shs(T) = \shs(\theta(T))$ for a binary total partition $T$.
Accordingly, we may formulate a recurrence for the numbers $h(n,k)$. Suppose that $U =
\tilde\gamma(T,X) \in \BT{[2,n]}$ (so $T \in \BT{[2,n-1]}$ has $2n-5$ nodes, of which $n-3$
are firstborn) and $|\shs(U)| = k$.  Then one of the following is true: either
    \begin{itemize}
    \item $|\shs(T)| = k$ and $X = \rt{T}$;
    \item $|\shs(T)| = k-1$ and $X$ is either one of the $n-3$ secondborn nodes of $T$, or
        one of the $k-1$ members of $\shs(T)$; or
    \item $|\shs(T)| = k-2$ and $X$ is one of the $(n-3)-(k-2)$ firstborn nodes
        of $T$ which do not belong to $\shs(T)$.
    \end{itemize}
Therefore $h(n,k)$ is defined by the recurrence
    \begin{equation} \label{recurhnk}
    h(n,k) ~=~ h(n-1,k) + (n+k-4) \: h(n-1,k-1) + (n-k-1) \: h(n-1,k-2)
    \end{equation}
with base cases $h(n,k) = 0$ if $k<0$ or $k>n-2$, $h(n,0) = 1$ for $n \geq 2$.

Let $n \in \Nn$.  A \defterm{matching} on $[1,2n]$ is a partition of $[1,2n]$ into $n$ pairs.  
A pair $\{i,j\}$ is \defterm{short} if $|i-j|=1$; otherwise it is \defterm{long}.
Define
    \begin{eqnarray*}
    \Match(n)   &=& \left\{ \text{matchings on}~ [1,2n] \right\}, \\
    \Match(n,k) &=& \left\{ X \in \Match(n) \st X ~\text{has}~ k ~\text{long pairs} \right\}, \\
    m(n,k)      &=& |\Match(n,k)|.
    \end{eqnarray*}
Kreweras and Poupard~\cite{KP74} gave recurrences and a closed formula for $m(n,k)$.  We will
use a slightly different argument
to show that $m(n,k)$ is given by a recurrence equivalent to \eqref{recurhnk}.

Note first that $m(n,0)=1$ for all $n$, because the only matching on $[1,2n]$ with no long pairs
is $\bs{ \{1,2\}, \{3,4\}, \dots, \{2n-1,2n\} }$.

Let $X \in \Match(n-1)$ and $2 \leq p \leq 2n$.  We define
a matching $X^p \in \Match(n)$ by \emph{inserting}
the pair $\{1,p\}$ into $X$ as follows.  First, we relabel $X$ to obtain a matching on $[1,2n]
\sm \{1,p\}$; this amounts to replacing each $x \in [1,2n-2]$ by $x+1$ if $x<p-1$, or by $x+2$
if $x \geq p-1$.  Having done this, we obtain $X^p$ by adjoining the pair $\{1,p\}$.

The map sending $(X,p)$ to $X^p$ is a bijection from $\Match(n-1) \x [2,2n] \to \Match(n)$, so by
induction $|\Match(n)| = (2n)!/2^n n!$. Moreover, we can derive a recurrence for $m(n,k)$.  Let $X
\in \Match(n,k)$.  If $p=2$, then $X^p \in \Match(n+1,k)$.  If $\{p-2,p-1\} \in X$
(which occurs for
$n-k$ values of $p$), then $X^p \in \Match(n+1,k+2)$ (because $\{p-2,p-1\}$ becomes the long pair
$\{p-1,p+1\}$ in $X^p$).  For the other $n+k$ values of $p$, we have $X^p \in \Match(n+1,k+1)$.
Accordingly, $m(n,k)$ is defined by the recurrence
    \begin{equation} \label{recurmnk}
    m(n,k) ~=~ m(n-1,k) + (n+k-2) \, m(n-1,k-1) + (n-k+1) \, m(n-1,k-2)
    \end{equation}
with base cases $m(n,k) = 0$ if $k<0$ or $k>n$, and $m(n,0) = 1$ for all $n\geq 0$.

\begin{theorem} \label{hvector-thm}
Let $n \geq 2$ and $0 \leq k \leq n-2$.  Then $h(n,k) = m(n-2,k)$, the number of matchings on
$[1,2n-4]$ with $k$ long pairs.
\end{theorem}

\begin{proof}
When $k=0$ or $n=2$, equality is immediate from the base cases.
Otherwise, assume that $h(n',k) = m(n'-2,k)$ for $n' < n$.  Then
\eqref{recurhnk} and \eqref{recurmnk} give
    \begin{equation*}
    \begin{split}
    m(n&-2,k) \\
    &=~ m(n-3,k) + (n+k-4) \: m(n-3,k-1) + (n-k-1) \: m(n-3,k-2) \\
    &=~ h(n-1,k) + (n+k-4) \: h(n-1,k-1) + (n-k-1) \: h(n-1,k-2) \\
    &=~ h(n,k).
    \end{split}
    \end{equation*}
\end{proof}


\section{A Recursive Lower Bound for Degree} \label{degree-section}

In this section, we return to geometry.  To complete the proof of
the main theorem, we must show that $J_n$ is the initial ideal of
an ideal defining the slope variety $\SVA(K_n)$ scheme-theoretically.
The key ingredient is to bound the (geometric) degree of $\SVA(K_n)$ from
below.  We do this by studying a family of \emph{flattened slope varieties}
$\SVA(n,k)$ forming a filtration of $\SVA(K_n)$: that is,
    $$\SVA(K_n) ~=~ \SVA(n,1) \supset \SVA(n,2) \supset \dots \supset \SVA(n,n).$$
We obtain a recursive lower bound $e(n,k) \leq \deg \SVA(n,k)$.  In Section~\ref{dpt-section},
we interpret the numbers $e(n,k)$ combinatorially and show that they actually
give the degree exactly.

We begin by recalling some standard facts about the geometric notion of degree.
Let $X \subset \Aa^N$ (or $\Pp^N$)
be an algebraic set of dimension $d$. The \defterm{degree} of $X$, denoted $\deg X$, is the
number of intersection points of $X$ with a generic affine linear subspace of codimension $d$.  
We extend this definition to locally closed sets $X,Y$ by putting $\deg X = \deg \overline{X}$.
Degree is multiplicative on products: $\deg (X \x Y) = (\deg X)(\deg Y)$.
Moreover, if $H$ is a hyperplane, then
$\deg (X \cap H) \leq \deg X$.  Finally,
if $X$ is the union of locally closed sets $C_1,\dots,C_k$
(for instance, if these are the irreducible components of $X$),
then the degree of $X$ is the sum of the degrees of all those $C_i$ whose
dimension equals $\dim X$.

We will also need some general facts about graph varieties; for details, see~\cite{JLM1}.
An \defterm{(affine) picture} $\Pic$ of $K_n$ consists of $n$ points $\Pic(1),\dots,\Pic(n)$
in the plane $\Aa^2$, and $\binom{n}{2}$ non-vertical lines $\Pic(1,2),\dots,\Pic(n-1,n)$,
subject to the conditions $\Pic(i),\Pic(j) \in \Pic(ij)$ for all $i,j$.  The set of all pictures is
the \defterm{affine picture variety} $\PVA(K_n)$.  
The affine slope variety $\SVA(K_n)$
is thus obtained by projecting $\PVA(K_n)$ on affine
coordinates
corresponding to the slopes of the lines $\Pic(ij)$.  We write $\phi$ for the
surjection $\PVA(K_n) \surj \SVA(K_n)$, and refer to a point in $\SVA(K_n)$
as a \defterm{slope picture} $\mm = (m_{ij})$.

A \defterm{partition} of a set $B$ is a collection of sets $\A = \{A_i\}$ such that $B$ is the
disjoint union of the $A_i$.  For each partition $\A$ of $[1,n]$, there is a corresponding
\defterm{cellule} $\PVA_{\A}(K_n)$ defined by
    \begin{equation*}
    \PVA_{\A}(K_n) \;=\; \left\{ \Pic \in \PVA(K_n) ~\left\vert~
    \begin{array}{l} \Pic(i)=\Pic(j) \text{ if and only if } \\
    i,j \text{ belong to the same part of } \A
    \end{array} \right. \right\}.
    \end{equation*}
Thus $\PVA(K_n)$ is the disjoint union of the cellules $\PVA_{\A}(K_n)$ as $\A$ ranges over 
all partitions.  Of particular importance are the \defterm{discrete cellule}, corresponding
to the partition of $[1,n]$ into $n$ singleton sets, and the \defterm{indiscrete cellule},
corresponding to the partition with only one part.

\begin{lemma} \label{trash-indiscrete-lemma}
Let $\mm \in \SVA(K_n)$.  Then there exists a picture $\Pic \in \phi^{-1}(\mm)$ in which not
all points $\Pic(1), \dots, \Pic(n)$ are the same.
Equivalently, $\phi$ remains surjective if the indiscrete cellule is deleted from
its domain.
\end{lemma}

\begin{proof}
The equations defining a picture of $G$ may be written in matrix form as $MX=0$, where $M$ is
a matrix whose entries are linear forms in the variables $m_e$ (as in
\eqref{deftreepolymatrix}) and $X$ is the column vector $[x_i-x_1]_{2 \leq i \leq n}$.  Every
maximal minor of $M$ corresponds to a subset of $E(K_n)$ of cardinality $2n-2$, which must
contain at least one rigidity circuit.  If $\mm \in \SVA(G)$, then every maximal minor of
$M(\mm)$ vanishes.  Therefore $M(\mm)$ has a nonzero nullvector $X$, which gives the
$x$-coordinates (up to translation) of a picture $\Pic$ not in the indiscrete cellule.
\end{proof}

For $1 \leq k \leq n$, define an algebraic subset of $\SVA(K_n)$ by
    \begin{equation}
    \SVA(n,k) = \{\mm \in \SVA(K_n) \st m_{ij} = 0 ~\text{for}~ 1 \leq i 
    < j \leq k\}.
    \end{equation}
Note that $\SVA(n,1) = \SVA(K_n)$.

\begin{proposition} \label{dim-slopevar-prop}
Let $n \geq 2$ and $1 \leq k \leq n$.  Then
    $$\dim\, \SVA(n,k) \: = \: \begin{cases}
      0 & ~\text{if } n=k, \\
      2n-k-2 & ~\text{if } n>k.
    \end{cases}$$
\end{proposition}

\begin{proof}
If $n=k$, then $\SVA(n,k)$ is a point.  If $k=1$, then $\dim \SVA(n,k) = \dim \SVA(K_n) = 2n-3 =
2n-k-2$.  These two observations complete the proof for $n=2$.  For $n>2$, we assume inductively
that the formula holds for all smaller $n$ and all $k$.  For each partition $\A$ of $[1,n]$, 
define a closed subset $F_{\A}$ of $\PVA_{\A}(K_n)$ by
    $$F_{\A} = \phi^{-1}(\SVA(n,k)) \,\cap\, \PVA_{\A}(K_n).$$
Note that Lemma~\ref{trash-indiscrete-lemma} implies that $\SVA(n,k) = \bigcup_{\A} 
\phi(F_{\A})$, where $\A$ ranges over all partitions with at least two parts.
In general, we may write $\A$ in the form
    $$\A ~=~ \left\{A_1 \dju B_1, ~ \dots, ~ A_s \dju B_s \right\}$$
where the $A_i$ (resp.\ $B_i$) form a partition of $[1,k]$ (resp.\ $[k+1,n]$), with
some $A_i$ (resp.\ $B_i$) allowed to be empty.  Let $a_i=|A_i|$ and $b_i=|B_i|$.  We may
order the parts of $\A$ in such a way that
    $$\begin{array}{llll}
    a_i > 0, & b_i > 0 & \text{for} & i \in [1,p]; \\
    a_i > 0, & b_i = 0 & \text{for} & i \in [p+1,q]; \\
    a_i = 0, & b_i = 1 & \text{for} & i \in [q+1,r]; \\
    a_i = 0, & b_i > 1 & \text{for} & i \in [r+1,s].
    \end{array}$$
for some $p \leq q \leq r \leq s$.  Note also that $q \leq k$.
A picture $\Pic \in F_{\A}$ is therefore given by $s$ distinct points in $\Aa^2$, of which $q$
lie on a common horizontal line, together with the slopes of the lines $\Pic(e)$ such that the
endpoints of the edge $e$ are in the same part of $\A$ and are not both $\leq k$.
This amounts to specifying $s$ $x$-coordinates and $s-q+1$ $y$-coordinates for the points
$\Pic(i)$; a point in $\SVA(a_i+b_i,a_i)$ for each $i \in [1,p]$; and a point in 
$\SVA(K_{b_i})$ for each $i \in [r+1,s]$. By the inductive hypothesis, we have
    \begin{eqnarray*}
    \dim \, F_{\A} &=& (2s-q+1)+ \sum_{i=1}^p (2b_i-a_i-2) + \sum_{i=r+1}^s (2b_i-3) \\
    &=& - 2p - s - q + 3r + 1 + 2\left(\sum_{i=1}^p b_i + \sum_{i=r+1}^s 
        b_i \right) - \sum_{i=1}^p a_i \\
    &=& 3r - 2p - q  - s + 1 + 2(n-k-(r-q)) - \sum_{i=1}^p a_i \\
    &=& 2n - 2k - 2p + q + r - s + 1 - \sum_{i=1}^p a_i \\
    &\leq& 2n - 2k - 2p + q + r - s + 1 \\
    &=& (2n-k+1) + (-k - 2p + q + r - s) \\
    &\leq& 2n-k+1
    \end{eqnarray*}
because $r \leq s$ and $q \leq k$.
Every fiber of $\phi$ has dimension at least $3$, so $\dim \phi(F_\A) \leq 2n-k-2$.
Since $\SVA(n,k)$ is the disjoint union of the $F_{\A}$, its dimension is at
most $2n-k-2$ as well.  On the other hand, if $\A$ is the discrete partition into
$n$ singleton sets, then $p=0$, $q=k$, and $r=s=n$.  So equality holds throughout the
preceding calculation, and $\dim \SVA(n,k) = 2n-k-2$.
\end{proof}

Define
    $$\SVAP(n,k) = \left\{ \mm \in \SVA(n,k) \st m_{1,k+1}=0 \right\}.$$
This set is the intersection of $\SVA(n,k)$ with a hyperplane, so its degree is a lower
bound for that of $\SVA(n,k)$.  We will calculate the degree of $\SVAP(n,k)$ by identifying its
irreducible components of maximal dimension.  First, we calculate the dimension of $\SVAP(n,k)$.

\begin{proposition} \label{another-slopevar-dim-prop}
Suppose $n>k \geq 1$ and $n \geq 2$.  Then $\dim \SVAP(n,k) = 2n-k-3$.
\end{proposition}

\begin{proof}
Intersecting $\SVA(n,k)$ with a codimension-$1$ hyperplane can lower the dimension by at most
$1$.  Hence for all $n$ and $k$, we have $\dim \SVAP(n,k) \geq 2n-k-3$.

For the reverse inequality, we induct on $n$. If $n=2$ and $k=1$, then $2n-k-3=0$, and indeed
$\SVAP(2,1) = \SVA(2,2)$ is a point.  Now suppose that $n>2$ and that the formula holds for
all smaller $n$ and all $k$.  For each partition $\A$ of $[1,n]$,
define a closed subset $F'_{\A}$ of $\PVA_{\A}(K_n)$ by
    $$F'_{\A} = \phi^{-1}\left(\SVAP(n,k)\right) \cap \PVA_{\A}(K_n).$$
By Lemma~\ref{trash-indiscrete-lemma}, $\SVAP(n,k) = \bigcup_{\A} \phi(F'_{\A})$,
where $\A$ ranges over all partitions of $n$ with at least two parts.  We need to show that
$\dim F'_{\A} \leq 2n-k$ for all such $\A$.  In general, we may write
    $$\A ~=~ \left\{ \{k+1\} \dju A_1 \dju B_1, ~ A_2 \dju B_2, ~ \dots, ~
      A_s \dju B_s \right\}$$
where the $A_i$ (resp.\ $B_i$) form a partition of $[1,k]$ (resp.\ $[k+2,n]$),
with some $A_i$ (resp.\ $B_i$) allowed to be empty.  Let $a_i=|A_i|$ and $b_i=|B_i|$.
Order the parts of $\A$ so that for some $p \leq q \leq r \leq s$,
    $$\begin{array}{llll}
    a_i > 0, & b_i > 0 & \text{for} & i \in [2,p]; \\
    a_i > 0, & b_i = 0 & \text{for} & i \in [p+1,q]; \\
    a_i = 0, & b_i = 1 & \text{for} & i \in [q+1,r]; \\
    a_i = 0, & b_i > 1 & \text{for} & i \in [r+1,s]. \\
    \end{array}$$
If $1 \in A_1$, then a picture $\Pic \in F'_{\A}$ is given by the following data: $s$ distinct 
points in $\Aa^2$, of which $q$ lie on a common horizontal line (so there are $s-q+1$ $x$-coordinates
and $s$ $y$-coordinates); a point in
$\SVAP(1+a_1+b_1,a_1)$; a point in $\SVA(a_i+b_i,a_i)$ for each $i \in [2,p]$;
and a point in $\SVA(K_{b_i})$ for each $i \in [r+1,s]$.  By induction,
    $$\dim\, F'_{\A} ~=~ (2s-q+1) + (a_1+2b_1-1) + \sum_{i=2}^p (a_i+2b_i-2)
        + \sum_{i=r+1}^s (2b_i-3).$$
Now a simple calculation, which we omit because of its similarity to that of
Proposition~\ref{dim-slopevar-prop}, shows that
$\dim\, F'_{\A} \leq 2n-k$ as desired.


If $1 \not\in A_1$, then a picture $\Pic \in F'_{\A}$ is given by the following data: $s$
distinct points in $\Aa^2$, of which $q$ lie on a common horizontal line; a point in
$\SVA(a_i+b_i,a_i)$ for each $i \in [1,p]$; and a point in $\SVA(K_{b_i})$ for each $i \in
[r+1,s]$.  Therefore
    \begin{eqnarray*}
    \dim\, F'_{\A} &=& (2s-q+1) + \sum_{i=1}^p (a_i+2b_i-2) + \sum_{i=r+1}^s (2b_i-3) \\
    &=& (-2p-q+3r-s+1) + \sum_{i=1}^p a_i
        + 2\left( \sum_{i=1}^p b_i + \sum_{i=r+1}^s b_i \right) \\
    &=& 2n-2k-2p+q+r-s-1 + \sum_{i=1}^p a_i.
    \end{eqnarray*}
Note that $\sum_{i=1}^p a_i \leq k-(q-p)$, since $\{A_1, \dots, A_q\}$ is a partition of
$[1,k]$ with no empty parts.  Therefore $\dim F'_{\A} ~\leq~ 2n-k-p+r-s-1 ~<~ 2n-k$ as
desired.
\end{proof}

\begin{theorem} \label{geom-recurrence-thm}
Let $n \geq 2$ and $k \in [1,n]$.  Then
    $$\deg \SVA(n,k) \geq e(n,k),$$
where $e(n,n)=1$ for all $n$; $e(2,1)=1$; and
    \begin{equation*}
    \begin{aligned}
    e(n,k)&~=~ e(n,k+1) ~+~ e(n-1,k-1) \\
    &+~ \sum_{t=1}^{k-1} \sum_{u=0}^{n-k-2} \binom{k-1}{t} \binom{n-k-1}{u} 
    ~ e(t+u+1,t) ~ e(n-t-u,k-t+1).
    \end{aligned}
    \end{equation*}
\end{theorem}

\begin{proof}
For all $n$, the variety $\SVA(n,n)$ is a point, so it has degree $1=e(n,n)$.  The variety
$\SVA(2,1)=\SVA(K_2)$ is a line, whose degree is again $1=e(2,1)$.  Now suppose that $n>k$ 
and $n>2$.  
As noted before, $\deg \SVAP(n,k) \leq \deg \SVA(n,k)$.
The idea of the proof is to write $\SVAP(n,k)$ as a finite union of locally closed
subsets.  Once we have done this, we may obtain a lower bound for $\deg \SVAP(n,k)$ by
summing lower bounds for the degrees of those locally closed
subsets of dimension $2n-k-3$.

First, note that $\SVAP(n,k) \supset \SVA(n,k+1)$, and that this second set has dimension
$2n-k-3$.

Second, let $\mm = (m_{ij})$ be a slope picture in $\SVAP(n,k) \sm \SVA(n,k+1)$.  
That is, $m_{ij}=0$ for $1 \leq i<j \leq k$, and $m_{1,k+1}=0$, but $m_{2,k+1}, \dots,
m_{k,k+1}$ are not all zero. Define two sets of vertices $T(\mm)$, $U(\mm)$ by
    \begin{eqnarray*}
    T(\mm) &=& \left\{ i \in [2,k] \st m_{i,k+1} \neq 0 \right\}, \\
    U(\mm) &=& \left\{ i \in [k+2,n] \st \Pic(i) = \Pic(k+1) ~\text{for all}~
    \Pic \in \phi^{-1}(\mm) \right\}.
    \end{eqnarray*}
Then $T(\mm) \neq \emptyset$, and
$\SVAP(n,k) \sm \SVA(n,k+1)$ is the disjoint union of the sets
    $$\STU = \left\{ \mm \in \SVAP(n,k) \sm \SVA(n,k+1) \st
    T(\mm)=T, ~ U(\mm)=U \right\}$$
for $\emptyset \neq T \subset [2,k]$ and $U \subset [k+2,n]$.  On $\phi^{-1}(\STU)$, one has
$\Pic(i)=\Pic(k+1)$ for all $i \in T \cup U$.  Thus the data for a picture in
$\phi^{-1}(\STU)$ consists of $n-|T|-|U|$ points, of which $k-|T|+1$ must lie on a common
horizontal line, together with the slopes of lines corresponding to edges with both endpoints
in $T \cup U$.  Applying the surjective map $\phi$, we see that $\STU$ has the form
    \begin{equation} \label{stu-prod}
    \STU \isom \SVA(n-t-u,k-t+1) \x \SVA(t+u+1,t)
    \end{equation}
where $t=|T|$ and $u=|U|$.  If $U \neq [k+2,n]$, then $n-t-u > k-t+1$.  Applying 
Proposition~\ref{dim-slopevar-prop} to \eqref{stu-prod}, one obtains
    $$\begin{array}{lll}
    \dim\, \STU &=& \left( 2(n-t-u)-(k-t+1)-2 \right) + \left( 2(t+u+1)-t-2 \right) \\
    &=& 2n-k-3.
    \end{array}$$
If $U = [k+2,n]$, then $n-t-u = k-t+1$.  Now the first factor in \eqref{stu-prod} is a 
single point, so applying Proposition~\ref{dim-slopevar-prop} yields
    $$\dim\, \STU ~=~ 2(t+u+1)-t-2 ~=~ (2n-k-3)+(t-k+1).$$
This is strictly less than $2n-k-3$ unless $T = [2,k]$, in which case $\STU \isom 
\SVA(n-1,k-1)$.  Putting all this together, we obtain
    \begin{align*}
    \deg\, \SVAP(n,k) \;\geq\; &\deg\, \SVA(n,k+1) + \deg\, \SVA(n-1,k-1) \\
    &+ \sum_{\substack{\emptyset \,\neq\, T \,\subset\, [2,k] \\
        U \,\subsetneq\, [k+2,n]}} \deg\, \SVA(n-t-u,k-t+1) ~ \deg\, \SVA(t+u+1,t).
    \end{align*}
Now, summing over $t$ and $u$ instead of $T$ and $U$, and multiplying the summand by the 
appropriate binomial coefficients to reflect the number of ways of choosing the sets $T$ and 
$U$, one obtains the desired recurrence.
\end{proof}


\section{Decreasing Planar Trees} \label{dpt-section}

The recurrence established in Theorem~\ref{geom-recurrence-thm} has a combinatorial
interpretation in terms of \defterm{decreasing planar trees}.  We begin with some general
facts about these objects, recalling the general definitions and terminology for rooted planar
trees given in Section~\ref{rpts-section}.

\begin{defn} \label{def-dpt}  
Let $V$ be a finite subset of $\Nn$.  A \defterm{decreasing planar tree on $V$} is a rooted
planar tree $T$ with the following property:  if $v,w$ are nodes of $T$ with $v$ an ancestor
of $w$, then $v>w$.  (In particular, $T$ is rooted at the maximum element of $V$.)
The set of all decreasing planar
trees on $V$ is denoted $\DPT(V)$; as usual we abbreviate $\DPT([1,n])$ by
$\DPT(n)$.
\end{defn}

The number of decreasing planar trees is given, once again, by the 
double factorial numbers:
$$
|\DPT(n)| ~=~ (2n-3)(2n-5) \cdots (3)(1) ~=~ 
\frac{(2n-2)!}{(n-1)!~2^{n-1}}.
$$
(This is an elementary combinatorial exercise; see, e.g.,~\cite[pp.~13--16]{Stanley-CCA}.)
We will obtain a recurrence enumerating decreasing planar trees on $V$ by the
\defterm{largest leaf} statistic:
$L(T) = \max \left\{ i \in V \st i~\text{is a leaf} \right\}$.
Let
    \begin{eqnarray*}
    \DPT(n,k) &=& \left\{ T \in \DPT(n) \st L(T) \leq k \right\}, \\
    \dpt(n,k) &=& |\DPT(n,k)|.
    \end{eqnarray*}
Note that $\dpt(n,1)=1$ for all $n$, since if $L(T)=1$ then the tree
$T$ can only be the path
in which each $j>1$ has the unique child $j-1$. Moreover, if $n \geq 2$, then
$\DPT(n,n) = \DPT(n,n-1) = \DPT(n)$.  In addition, $\dpt(n,k)=0$ if $n>0$
and $k \notin [1,n]$.  Conventionally, we put $\dpt(0,k)$ for all $k$.

If $T_1$ and $T_2$ are rooted planar trees and $X$ is a node of $T_2$, we may graft $T_1$
to $T_2$ by attaching it as the oldest subtree of $X$.  We denote the resulting tree by
$T_1 \graftat{X} T_2$,
suppressing the subscript if $X=\rt{T_2}$ (the most common case).
To illustrate this operation, let $T$ be the binary tree
pictured in Example~\ref{decomp-tree-example}.
Let $T_1$ be the subtree rooted
at the node labeled $245$, and let $T_2$ be the tree obtained from $T$ by
deleting $T_1$.  Then $T = T_1 \graftat{X} T_2$, where $X$ is the node
labeled $2456$.

Suppose that $T_1 \in \DPT(V_1)$ and $T_2 \in \DPT(V_2)$, where $V_1 \cap V_2 =
\emptyset$.  If $X \in V_2$ and $\rt{T_1} < X$,
then $T_1 \graftat{X} T_2 \in \DPT(V_1 \cup V_2)$.  In 
addition, every tree $T$ with more than one node can be written uniquely as
$T_1 \graft T_2$, by taking $T_1 = T|{\rt{T}}^{(1)}$ and $T_2=T-T_1$.

\begin{lemma} \label{dpt-lemma-one}
Let $k \in [2,n]$.  Then
    $$\dpt(n,k-1) ~=~ \sum_{a=1}^{k-1} \binom{k-1}{a} ~ \dpt(a,a) ~ \dpt(n-a,k-a).$$
\end{lemma}

\begin{proof}
Let $T_1$ be a decreasing planar tree with nodes $A \subset
[1,k-1]$, $|A|=a<k$.  There are $\binom{k-1}{a}\dpt(a,a)$ such trees.  Meanwhile, let $T_2$
be a decreasing planar tree with nodes $[a+1,n]$ and largest leaf $L(T_2) \leq k$.
Subtracting $a$ from the label of every node gives a bijection between such trees $T_2$
and the set $\DPT(n-a,k-a)$.  Therefore the number of pairs $(T_1,T_2)$
is counted by the right-hand side of the desired equality.  For each such pair, let
    $$f(T_1,T_2) = T_1 \graftat{k} T_2.$$
The map $f$ is clearly injective, and $f(T_1,T_2) \in \DPT(n,k-1)$.  On the other hand, for
any $T \in \DPT(n,k-1)$, let $T' = T|k^{(1)}$ and $T''=T-T'$.
Then $T' \graftat{k} T'' = T$.  So $f$ maps onto
$\DPT(n,k-1)$, and the desired equality follows.
\end{proof}

This formula may be taken as a recursive definition for the numbers $\dpt(n,k)$.
For small values of $n$ and $k$, they are given by the following table.
$$\begin{array}{c|cccccc}
    & k=1 & 2 & 3 & 4 & 5 & 6 \\ \hline
n=2 & 1 & 1\\
3   & 1 & 3 & 3\\
4   & 1 & 7 & 15 & 15\\
5   & 1 & 15 & 57 & 105 & 105\\
6   & 1 & 31 & 105 & 561 & 945 & 945
\end{array}$$

\vskip 0.1in
\begin{lemma} \label{dpt-lemma-two}
Let $1 < k \leq n$.  Then
    \begin{multline*}
    \sum_{a=1}^{k-1} \sum_{c=1}^{n-k-1} \binom{k-1}{a} \binom{n-k-1}{c}
    ~ \dpt(a+c,a) ~ \dpt(n-a-c,k-a) ~= \\
    \sum_{w=0}^{k-2} \sum_{y=1}^{n-k-1} \binom{k-1}{w} \binom{n-k-1}{y}
    ~ \dpt(w+y,w+1) ~ \dpt(n-w-y,k-w-1).
    \end{multline*}
\end{lemma}

\begin{proof} 
First, we interpret the sums on either side of the desired equality
as enumerating certain kinds of decreasing planar trees.
The left-hand side counts trees of the form $T =
T' \graft T'' \in \DPT(n)$ such that
    \begin{equation} \label{typetwo}
    \begin{aligned}
    V(T')  &= A \dju C,              \qquad & L(T')  &\in A, \\
    V(T'') &= B \dju D \dju \{k,n\}, \qquad & L(T'') &\in B \dju \{k\},
    \end{aligned}
    \end{equation}
where
    \begin{equation*}
    \begin{aligned}
    A \dju B &= [1,k-1],   \qquad & 1 \;\leq\; a &=|A| \;\leq\; k-1, \\
    C \dju D &= [k+1,n-1], \qquad & 1 \;\leq\; c &=|C| \;\leq\; n-k-1.
    \end{aligned}
    \end{equation*}
Meanwhile, the right-hand side counts trees of the form $U = U' \graft U'' \in \DPT(n)$
such that
    \begin{equation} \label{typethree}
    \begin{aligned}
    V(U')  &= W \dju Y,              \qquad & L(T')  &\in W \dju \{\min(Y)\}, \\
    V(U'') &= X \dju Z \dju \{k,n\}, \qquad & L(T'') &\in X,
    \end{aligned}
    \end{equation}
where
    \begin{equation*}
    \begin{aligned}
    W \dju X &= [1,k-1],   \qquad & 0 \leq w&=|W| \;\leq\; k-2, \\
    Y \dju Z &= [k+1,n-1], \qquad & 1 \leq y&=|Y| \;\leq\; n-k-1.
    \end{aligned}
    \end{equation*}

Second, we show how to perform ``surgery'' on a tree of the
form~\eqref{typetwo} to obtain a tree of the form~\eqref{typethree}.
Suppose that $T = T' \graft T''$ satisfies \eqref{typetwo}.  Define a tree $f(T)$ as follows.

\noindent $\bullet$ \quad {\it If $k$ has children in $T$\/}, then $f(T)=T$.  In this case $B
\neq \emptyset$, so $A \subsetneq [1,k-1]$.  Also, $L(T) < k$; in particular $L(T'') \in B$.  
Therefore, $T$ satisfies \eqref{typethree}, with $W=A$, $X=B$, $Y=C$, $Z=D$.

\noindent $\bullet$ \quad {\it If $k$ is a leaf\/}, then we form $f(T)$ by detaching all
subtrees of $\min(C)$ and reattaching them to $k$ in the same birth order.
Note that there is at least one such subtree since $\min(C) > L(T)$.
Also, $L(f(T)) = \min(C)$, and $f(T)$
satisfies \eqref{typethree}, with
    \begin{equation*} 
    \begin{aligned}
    W &= A \sm  \{\text{descendants of}~ \min(C)\}, \qquad & Y &= C, \\
    X &= B \dju \{\text{descendants of}~ \min(C)\}, \qquad & Z &= D.
    \end{aligned}
    \end{equation*}

Third, we show how this surgery can be reversed.
Suppose that $U = U' \graft U''$ satisfies \eqref{typethree}.  Define a tree
$g(U)$ as follows:

\noindent $\bullet$ \quad {\it If $\min(Y)$ has children in $U$\/}, then $g(U)=U$.  In this
case $W \neq \emptyset$.  Also, $L(T) < k$; in particular $L(T'') \in X$.  Therefore, $T$
satisfies \eqref{typethree}, with $A=W$, $B=X$, $C=Y$, $D=Z$.

\noindent $\bullet$ \quad {\it If $\min(Y)$ is a leaf\/}, then we form $g(U)$ by detaching
all subtrees of $k$ and reattaching them to $\min(Y)$ in the same birth order.  Note that
since $k > L(U)$, this step is not trivial.  Also, $L(g(U)) = k$, and $g(U)$ satisfies
\eqref{typetwo}, with
    \begin{equation*} 
    \begin{aligned}
    A &= W \dju \{\text{descendants of}~ \min(Y)\}, \qquad & C &= Y, \\
    B &= X \sm  \{\text{descendants of}~ \min(Y)\}, \qquad & D &= Z.
    \end{aligned}
    \end{equation*}

Finally, we show that the functions $f$ and $g$ are mutual inverses.  Suppose that $T$ 
satisfies
\eqref{typetwo}.  In particular, $\min(C)$ has children in $T$, so if $f(T)=T$, then
$g(f(T))=T$.  On the other hand, if $k$ is a leaf in $T$, then by construction
$\min(C)=\min(Y)$ is a leaf in $f(T)$, and the two surgeries described above are
inverses by definition.  So $g(f(T))=T$ in any case.  The proof that $f(g(U))=U$ is
analogous.  So $f$ and $g$ are bijections and we are done.
\end{proof}

\begin{lemma} \label{dpt-lemma-three}
Let $1 < k \leq n$.  Then
    \begin{equation*}
    \dpt(n,k) = \dpt(n-1,k) +
        \sum_{a=1}^{k-1} \sum_{c=0}^{n-k-1} \binom{k}{a} \binom{n-k-1}{c}
        \dpt(a+c,a) ~ \dpt(n-a-c,k-a).
    \end{equation*}
\end{lemma}

\begin{proof}
Deleting the node $n$ gives a bijection between trees $T \in \DPT(n,k)$
such that $n$ has exactly one child, and $\DPT(n-1,k)$.  This accounts
for the term $\dpt(n-1,k)$ on the right-hand side.
On the other hand, suppose that
$T \in \DPT(n,k)$ is a
tree in which $n$ has more than one child.  Then
$T = T' \graft T''$, where
    \begin{equation*}
    \begin{aligned}
    V(T')    &= A \dju C,            \qquad & L(T')    &\in A,\\
    V(T'')   &= B \dju D \dju \{n\}, \qquad & L(T'')   &\in B,\\
    A \dju B &= [1,k],               \qquad & 1 \leq a &= |A| \leq k-1, \\
    C \dju D &= [k+1,n-1],           \qquad & 0 \leq c &= |C| \leq n-k-1.
    \end{aligned}
    \end{equation*}
Conversely, if a pair $(T',T'')$ satisfies these conditions, then $T = T' \graft
T''$ belongs to $\DPT(n,k)$, and $n$ has at least two children in $T$.  These pairs are
enumerated by the double sum on the right-hand side of the desired equality.
\end{proof}

Applying Pascal's identity
    $$\binom{k}{a} = \binom{k-1}{a} + \binom{k-1}{a-1}$$
to the equation of Lemma~\ref{dpt-lemma-three}, we obtain
    \begin{equation} \label{d-rec-two}
    \begin{split}
    &\dpt(n,k) ~=~ \dpt(n-1,k) \\
    &+ \sum_{a=1}^{k-1} \sum_{c=0}^{n-k-1} \binom{k-1}{a} \binom{n-k-1}{c}
    ~ \dpt(a+c,a) ~ \dpt(n-a-c,k-a) \\
    &+ \sum_{a=1}^{k-1} \sum_{c=0}^{n-k-1} \binom{k-1}{a-1} \binom{n-k-1}{c}
    ~ \dpt(a+c,a) ~ \dpt(n-a-c,k-a).
    \end{split}
    \end{equation}

Breaking off the $c=0$ term of the first double sum in \eqref{d-rec-two} yields
    \begin{equation} \label{d-rec-three}
    \begin{split}
    &\dpt(n,k) ~=~ \dpt(n-1,k) 
    ~+~ \sum_{a=1}^{k-1} \binom{k-1}{a} ~ \dpt(a,a) ~ \dpt(n-a,k-a) \\
    &+ \sum_{a=1}^{k-1} \sum_{c=1}^{n-k-1} \binom{k-1}{a} \binom{n-k-1}{c}
    ~ \dpt(a+c,a) ~ \dpt(n-a-c,k-a) \\
    &+ \sum_{a=1}^{k-1} \sum_{c=0}^{n-k-1} \binom{k-1}{a-1} \binom{n-k-1}{c}
    ~ \dpt(a+c,a) ~ \dpt(n-a-c,k-a).
    \end{split}
    \end{equation}

Applying Lemma~\ref{dpt-lemma-one} to the single sum and Lemma~\ref{dpt-lemma-two} to 
the first double sum in \eqref{d-rec-three}, we obtain
    \begin{equation} \label{d-rec-four}
    \begin{split}
    &\dpt(n,k) ~=~ \dpt(n-1,k) + \dpt(n,k-1) \\
    &+ \sum_{w=0}^{k-2} \sum_{y=1}^{n-k-1} \binom{k-1}{w} \binom{n-k-1}{y}
    ~ \dpt(w+y,w+1) ~ \dpt(n-w-y,k-w-1) \\
    &+ \sum_{a=1}^{k-1} \sum_{c=0}^{n-k-1} \binom{k-1}{a-1} \binom{n-k-1}{c}
    ~ \dpt(a+c,a) ~ \dpt(n-a-c,k-a).
    \end{split}
    \end{equation}

In the first double sum of \eqref{d-rec-four}, we may change the upper limit on
$y$ from $n-k-1$ to $n-k$, because the additional summand is zero; in addition put
$x = n-k-y$. In the second double sum, put $w=a-1$ and $x=n-k-1-c$.  We obtain
    \begin{equation} \label{d-rec-five}
    \begin{split}
    &\dpt(n,k) ~=~ \dpt(n-1,k) \:+\: \dpt(n,k-1) \\
    &+ \sum_{w=0}^{k-2} \sum_{x=0}^{n-k-1} \binom{k-1}{w} \binom{n-k-1}{x-1}
    \, \dpt(n-k+w-x,w+1) \, \dpt(k-w+x,k-w-1) \\
    &+ \sum_{w=0}^{k-2} \sum_{x=0}^{n-k-1} \binom{k-1}{w} \binom{n-k-1}{x}
    \, \dpt(n-k+w-x,w+1) \, \dpt(k-w+x,k-w-1)
    \end{split}
    \end{equation}

Now, another application of Pascal's identity, together with the base cases
mentioned previously, yields the following:

\begin{proposition} \label{dpt-recurrence-prop}
The numbers $d(n,k)$ satisfy the recurrence
    \begin{multline*}
    \dpt(n,k) ~=~ \dpt(n-1,k) ~+~ \dpt(n,k-1) \\
    +~ \sum_{w=0}^{k-2} \sum_{x=0}^{n-k-1} \binom{k-1}{w} \binom{n-k}{x}
    ~\dpt(n-k+w-x,w+1) \; \dpt(k-w+x,k-w-1).
    \end{multline*}
\end{proposition}

Finally, we come to the geometric reason for all this enumeration.
The recurrence given by Proposition~\ref{dpt-recurrence-prop}
is easily turned into one giving the lower bounds for the degree
of the flattened slope variety $\SVA(n,k)$, as we now show.

\begin{theorem} \label{degree-bound-thm} We have
$\displaystyle \deg \SVA(n,k) \geq \dpt(n-1,n-k)$. In particular,
    $$\deg \SVA(K_n) \geq \dpt(n-1,n-1) = |\DPT(n-1)|.$$
\end{theorem}

\begin{proof}
Replacing $\dpt(i,j)$ with $e(i+1,i-j+1)$ in Proposition~\ref{dpt-recurrence-prop}, we obtain
    \begin{equation*}
    \begin{split}
    &e(n+1,n-k+1) ~=~ e(n,n-k) + e(n+1,n-k+2) \\
    & \qquad +~ \sum_{w=0}^{k-2} \sum_{x=0}^{n-k-1} \binom{k-1}{w} \binom{n-k}{x} \\
    & \qquad\qquad \x~ e(n-k+w-x+1,n-k-x) ~ e(k-w+x+1,x+2).
    \end{split}
    \end{equation*}
Now, setting
    $$n = n'-1, \qquad k = n'-k', \qquad x = k'-w'-1, \qquad w = x'$$
and removing the primes from the new variables, we obtain the recurrence of
Theorem~\ref{geom-recurrence-thm}.  It follows that $\dpt(n,k) = e(n+1,n-k+1)$ for all
$n$ and $k$, or equivalently $e(n,k) = \dpt(n-1,n-k)$.
\end{proof}


\section{Proof of the Main Theorem} \label{proof-main-section}

With all the pieces in place, we can now prove our main result.

\begin{proof}[Proof of Theorem~\ref{main-thm}]
The variety $\SVA(K_n)$ is irreducible and reduced, hence is defined scheme-theoretically by
the prime ideal $\sqrt{I_n}$.  In addition, we have
    \begin{equation} \label{inideals}
    J_n \subset \ini(I_n) \subset \ini(\sqrt{I_n}).
    \end{equation}
By Theorem~\ref{facet-decomp-thm} we have
    $$\codim J_n = 2n-3 = \dim \SVA(K_n) = \codim \sqrt{I_n}$$
and by Theorem~\ref{facet-count-thm} and Theorem~\ref{degree-bound-thm} we have
    $$\deg J_n = \frac{(2n-4)!}{2^{n-2}(n-2)!} \leq \deg \sqrt{I_n}.$$
We have seen in Theorem~\ref{shelling-thm} that the Stanley-Reisner complex of $J_n$ is
shellable. Hence $J_n$ is Cohen-Macaulay~\cite[Theorem~5.1.13]{BH}. In particular it is
unmixed, so any strictly larger ideal has either larger codimension or smaller degree.  
Meanwhile, $\ini(\sqrt{I_n})$ has the same codimension and degree as $\ini(\sqrt{I_n})$
by~\cite[Theorem~15.26]{Eisenbud}, Therefore equality holds in \eqref{inideals}. Since $J_n =
\ini(I_n)$ is Cohen-Macaulay, so is $I_n$~\cite[Corollary~B.2.1]{Vasconcelos}.  It follows that
$\SVA(K_n)$ is Cohen-Macaulay.  Finally, the formula for the Hilbert series is a
restatement of the results of Section~\ref{hvector-section}.
\end{proof}

We remark that for the degree bounds in Theorem~\ref{degree-bound-thm} to be sharp, equality
must hold in the lower bounds of Theorem~\ref{geom-recurrence-thm}.  That is, the flattened
slope variety $\SVA(n,k)$ has degree exactly $d(n-1,n-k)$.

We conclude by mentioning some problems for further study.
Using the computer algebra system {\it Macaulay\/}~\cite{Macaulay}, we have verified that
$I_n$ is generated by the (cubic) tree polynomials of the
$\binom{n}{4}$ copies of $K_4$ arising as subgraphs of $K_n$,
for $n\leq 9$.  Accordingly, we conjecture
that they do so for all $n$.  The difficulty is that there appears to be
no canonical way to write a tree polynomial $\TP(W)$ as an $R_{K_n}$-linear combination
of the $\TP(K_4)$'s, where $W$ is a wheel with five or more vertices.  Further computations
using {\it Macaulay\/} suggest that the wheel polynomials may form a \defterm{universal
Gr\"obner basis} for $I_n$ (that is, a Gr\"obner basis with respect to \emph{every} term
ordering), but it is not clear how to prove this.

The technique of using the Stanley-Reisner complex to obtain information about graph varieties
may be applicable to graphs other than $K_n$.  The likeliest candidates are graphs
with many symmetries,
such as threshold graphs or complete bipartite graphs.  The degree
and Hilbert series of the slope variety of a graph are combinatorial invariants
for which it would be nice to have a more elementary graph-theoretic description.
In addition, it would be of interest to determine whether $\SVA(G)$ is always defined
scheme-theoretically by the tree polynomials of rigidity circuit subgraphs
of $G$, and whether it is always Cohen-Macaulay.  (We have
computed $I_G$ for several small graphs; in all cases the Cohen-Macaulay property is
satisfied.)

An anonymous referee has pointed out that the space $\mathscr{T}_n$ of
\defterm{phylogenetic trees on $n$ vertices},
constructed by Billera, Holmes and Vogtmann in~\cite{BHV},
has degree $(2n-4)!/(2^{n-2}(n-2)!)$, the same as that of $\SVA(K_n)$.
Moreover, for $0<m<n$, the space $\mathscr{T}_n$ contains
many subspaces isomorphic to $\mathscr{T}_m$ \cite[p.~743]{BHV}, just as
$\SVA(K_n)$ contains many copies of $\SVA(K_m)$ (see Section~\ref{degree-section}).
It would be interesting to explore possible connections between pictures
and phylogenetic trees.

Other spaces related to graph varieties include the Fulton-Macpherson
\defterm{compactification of configuration space}~\cite{FM94} and the
De~Concini-Procesi \defterm{wonderful model of subspace arrangements}
\cite{DP95}.  The results of this article might serve as a starting point
for studying these relationships in more detail.
\bibliographystyle{plain}
\bibliography{biblio}
\end{document}